\documentclass[final]{siamltex}

\usepackage{amsmath}
\usepackage{amsfonts}
\usepackage{amssymb}
\usepackage{latexsym}
\usepackage{exscale,cmmib57}
\usepackage{graphicx,color}
\usepackage{multirow}
\usepackage{graphicx}
\usepackage{epic}
\usepackage{subfigure}
\usepackage{array}
\usepackage{epic,eepic}
\usepackage{fancybox}
\usepackage{psfrag}
\usepackage{caption}
\usepackage{float}
\usepackage{multirow}

\newcommand{\mean}[1]{\{ \kern -1.6mm \{#1\} \kern -1.6mm \}}
\newcommand{\jump}[1]{[ \kern -.7mm [#1] \kern -.7mm ]}

\newcommand{\be}{\begin{equation}}\newcommand{\ee}{\end{equation}}

\newcommand{\ndg}[1]{| \kern -.25mm \|{#1}| \kern -.25mm \|}

\renewcommand{\ldots}{\dotsc}

\title{Multilevel Sparse Kernel-Based Interpolation}

\author{Emmanuil H.
Georgoulis\thanks{Department of Mathematics, University of
Leicester, University Road, Leicester LE1 7RH, United Kingdom
(Emmanuil.Georgoulis@le.ac.uk).} \and Jeremy
Levesley\thanks{Department of Mathematics, University of Leicester,
University Road, Leicester LE1 7RH, United Kingdom
(j.levesley@le.ac.uk).} \and Fazli Subhan\thanks{Department of
Mathematics, University of Leicester, University Road, Leicester LE1
7RH, United Kingdom ((a):fs83@le.ac.uk,
(b):f.subhan@hotmail.com.)} }

\begin{document}

\maketitle

\hyphenation{di-men-sio-nal}

\begin{abstract}
A multilevel kernel-based interpolation method, 
suitable for moderately high-dim\-ensional function interpolation problems, is proposed. The method, termed \emph{multilevel sparse kernel-based interpolation} (MLSKI, for short), uses both level-wise and direction-wise multilevel decomposition of structured  (or mildly unstructured) interpolation data sites in conjunction with the
application of kernel-based interpolants with different scaling in
each direction. The multilevel interpolation algorithm is based on a hierarchical decomposition of the data sites, whereby at each level the detail is added to the interpolant by interpolating the resulting residual of the previous level. On each level, anisotropic radial basis functions are used for solving a number of small interpolation problems, which are subsequently linearly combined to produce the interpolant. MLSKI can be viewed as an extension
of $d$-boolean interpolation (which is closely related to ideas in sparse grid and hyperbolic crosses literature) to kernel-based functions, within the hierarchical multilevel framework to achieve accelerated convergence. Numerical experiments suggest that the new algorithm is numerically stable and efficient
for the reconstruction of large data in $\mathbb{R}^{d}\times \mathbb{R}$, for  $d = 2, 3, 4$,  with tens or even hundreds of thousands data points.  Also, MLSKI appears to be generally superior over classical radial basis function methods in terms of complexity, run time and convergence at least for large data sets. 
\end{abstract}



\pagestyle{myheadings} \thispagestyle{plain} \markboth{E.~H.~GEORGOULIS,  J.~LEVESLEY  AND F.~SUBHAN}{MULTILEVEL SPARSE KERNEL-BASED INTERPOLATION}

\section{Introduction}\label{sec1}

Over the last four decades, radial basis functions (RBFs) have been successfully applied to (scattered) data interpolation/approximation in $\mathbb{R}^{d}\times \mathbb{R}$ (see, e.g., \cite{Wend05} and the references therein for a literature review). The interest on kernel-based and, in particular, on RBF interpolants can be traced in their ability to produce global interpolants  of user-defined smoothness  without the shortcomings of multivariate polynomial interpolation. These interpolants admit generally  good convergence properties and they can be implemented in 
(essentially) dimension-independent fashion, making them potentially attractive for a number of applications.
 
Despite the above attractive properties, RBF interpolation can be cumbersome in practice. Solving the resulting linear system is challenging due to both the density and the ill-conditioning of the resulting interpolation matrix (see, e.g., \cite{Fasshauerbook2007,Schaback1995C}). A number of techniques have been proposed to deal with the ill-conditioning of the interpolation system. 

For smooth basis functions such as the Gaussian and the multiquadric, the ill-conditioning depends crucially on a width parameter.  Novel QR-based algorithms which eliminate the ill-conditioning problem for thousands of points have been developed  \cite{FB04,BEN, FasshauerandMccourt2011}. An alternative family of methods, where small local problems are solved which generate approximate cardinal functions has also been proposed \cite{Powell99,Beatson99}. The matrix associated with such an approach is much better conditioned, and allows for rapid solution using iterative methods. Finally, in \cite{MR2177889}, a stable interpolation algorithm on carefully selected nodes for Gaussian basis functions is given.

The introduction of RBFs with compact support \cite{wu,wendland} aims to address the density issue of the interpolation matrix. Moreover, a number of techniques have been developed to reduce the complexity of calculating the interpolant, involving multipole type expansion for a variety of RBFs  \cite{Beatson99}.   Using such methods is possible to compute an RBF interpolant with $O(k \log k) $-computations for quasi--uniform data, where $k$ is the number of data sites though, to the best of our knowledge, we are not aware of any methods which guarantee a bounded number of iterations independently of $k$. 
 
Thus, the complexity of kernel-based interpolation remains a challenge when $d\ge2$, as the number of data points is required to grow exponentially with respect to the dimension $d$ to ensure good convergence rates.  It is not surprising, therefore, that the use of RBFs in practice has been largely limited to tens of thousands of data sites, which, in turn, has restricted their application to low $d$, typically $d=1,2$ or $3$.

This work is concerned with the introduction of a kernel-based interpolation method on structured or mildly unstructured data sites, which aims to address the computational complexity issues of RBF interpolation for $d\ge 2$, while simultaneously reducing the ill-conditioning of the resulting interpolation problem. The new scheme, termed \emph{multilevel sparse kernel-based interpolation} (MLSKI), is based on a hierarchical decomposition of the data sites, whereby at each level the detail is added to the interpolant by interpolating the resulting residual of the previous level. On each level, anisotropic radial basis functions are used for solving a number of small interpolation problems, which are subsequently linearly combined to produce the interpolant; the new method can be viewed as an extension of (and, indeed, it has been inspired from)
the idea of $d$-boolean interpolation \cite{DelvosFJ1982, SS99, KW05b,HeglandandGarckeandChallis2007,GarckeandHegland2009} to kernel-based functions, which, in turn, is closely related to ideas in sparse grid \cite{Zen91,GSZ92, BGRZ94, HEM95,Gri98,GK03} and hyperbolic crosses \cite{Smo63,Beb60} literature. We note that
in~\cite{SchreiberAnja2000} hyperbolic cross products of one
dimensional RBFs have been considered. The hierarchical multilevel framework used to achieve accelerated convergence is relatively standard in the RBF literature \cite{Floater&Iske1996,Iske2001,IskeANDLeveslely2005,Hales&Levesley2002}.

In the simplest setting, the MLSKI algorithm assumes that the data sites are on a Cartesian uniform grid of size $N^d$ in $\mathbb{R}^d$, with $N=2^n+1$, for a final level $n\in\mathbb{N}$. For each level $0\le l\le n$ of the MLSKI algorithm, we construct a \emph{sparse kernel-based interpolant} of the interpolation residual as follows. We consider $O(l^{d-1})$ carefully chosen subsets of the data points of level $l$, each subset having size $O(2^l)$ data points. On each of these subsets, which we shall refer to as \emph{partial grids}, we solve the interpolation problem. As the partial grids are anisotropic in nature, we employ appropriate anisotropically scaled kernels (anisotropic RBFs \cite{Ani06,Ani07,Ani10}) for each interpolation problem on each partial grid. Once all the $O(l^{d-1})$ interpolants on the partial grids have been computed they are linearly combined to give the total interpolant (of the residual) for level $l$. Hence, the complexity of the resulting MLSKI algorithm is dominated by the complexity of the last step, i.e., one needs to solve $O(n^{d-1})$ interpolation problems of size $N$ and linearly combine the resulting interpolants on the partial grids. This, in turn, implies that the MLSKI algorithm admits (at least theoretically) $O(\sigma(N)\log^{d-1} \sigma(N))$-complexity, where $N$ is the number of grid points in \emph{each} direction, and $\sigma(N)$ is the complexity of the univariate $N$-point algorithm.  

The computation of each interpolant on the partial grids is completely independent from the other interpolants on each level $l$. Therefore, the MLSKI algorithm is ideally suited for implementation on parallel computers.  The MLSKI algorithm is tested in practice for a number of relevant test cases with $d=2,3,4$. The numerical experiments suggest that MLSKI is numerically stable in practice and efficient for the reconstruction of large data in $\mathbb{R}^{d}\times \mathbb{R}$, for  $d = 2, 3, 4$,  with hundreds of thousands of data points.  Also, MLSKI appears to be generally superior over classical radial basis function methods in terms of complexity, run time and convergence, at least for large data sets.

The remaining of this paper is organized as follows. In Section \ref{sec2},
we discuss anisotropic versions of radial basis functions, suitable for interpolation on data sites with anisotropic characteristics. In Section~\ref{sec3}, we introduce the sparse kernel-based interpolation method, which will be used in Section~\ref{sec4} at each level of the multilevel algorithm. The stability and implementation of the MLSKI algorithm is discussed in Section~\ref{sec5}. A series of numerical experiments is given in Section~\ref{sec6}, for $d$ = 2, 3, 4.  In
Section~\ref{sec7}, we draw some conclusions on the new developments
presented in this work and discuss possible extensions.

\section{Anisotropic RBF interpolation}\label{sec2}
\label{sec:Anisotropic radial basis functions interpolation}
Radial basis functions are radially symmetric by construction having hyper-spheres as their level surfaces. However, interpolation of data with anisotropic distribution of data sites in the domain requires special consideration. To this end, anisotropic
radial basis functions (ARBFs) have been introduced and used in practice~\cite{Ani06,Ani07,Ani10}; they are also known as  elliptic basis functions as they have hyper-ellipsoidal level surfaces.

\begin{definition}\label{Def:Anisotropic radial basis functions}
Let $\varphi(\|\cdot-\mathbf{x}_{j}\|)$ be a given RBF centred at $\mathbf{x}_{j}\in \mathbb{R}^{d}$ and  let $A\in\mathbb{R}^{d\times d}$ be an invertible matrix.
The anisotropic radial basis function
$\varphi_{A}$\label{varphiA} is defined by
$\varphi_{A}(\| \cdot-\mathbf{x}_{j}\|)=\varphi(\| A(\cdot-\mathbf{x}_{j})\|)$.
\end{definition}

\begin{figure}
\begin{center}
\includegraphics[width=10cm,height=5cm]{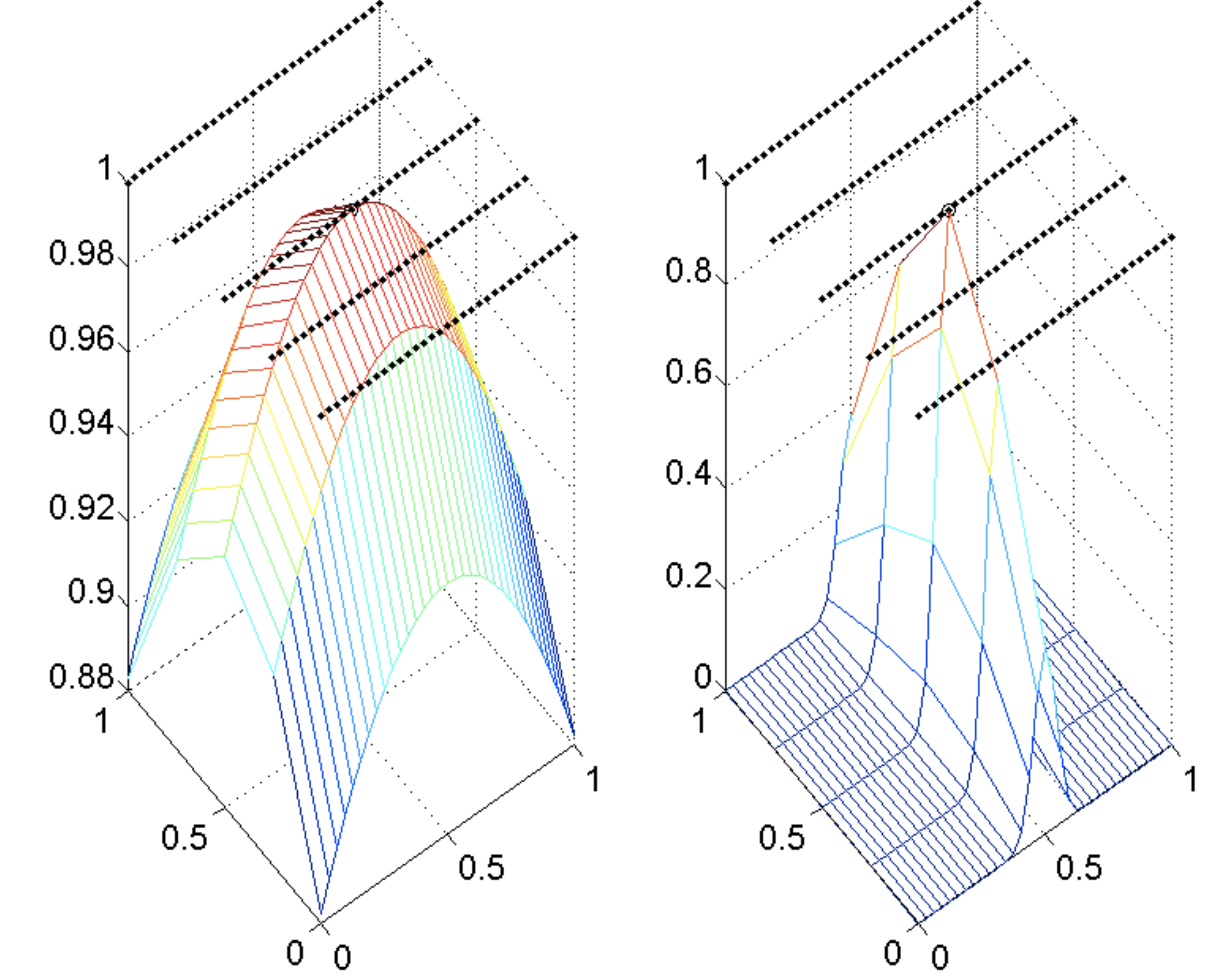}
\end{center}
\caption{Gaussian and anisotropic RBFs,  with shape parameter $c=0.5$.}
\label{Fig:GaussianAnisotropic}
\end{figure}

Evidently, $\varphi_{A}(\|\cdot\|)\equiv\varphi(\|\cdot\|)$  when $A$ is the $d\times d$-identity matrix.
In Figure~\ref{Fig:GaussianAnisotropic}, the Gaussian RBF and the corresponding
anisotropic Gaussian RBF centred at $(0.5,0.5)$ for  $A={\rm diag} (2^5, 2^2)$ are drawn. (Here, we use the notation ${\rm diag}(v)\in\mathbb{R}^{d\times d}$ for a diagonal matrix with diagonal entries are given by the components of the vector $v\in\mathbb{R}^d$.)

We now discuss the the solution of an anisotropic interpolation
problem. We restrict the discussion to positive definite kernels, which will be used later, though anisotropic versions of conditionally positive definite REBs is completely completely analogous.

For data sites 
$X:=\{\mathbf{x}_{1},\ldots,\mathbf{x}_{N}\}$ contained in a computational domain $\Omega\subset\mathbb{R}^{d}$, we consider the interpolation data $\{(\mathbf{x}_{i},y_{i}):\mathbf{x}_{i} \in
X,i=1,\ldots,N\}$. Let $\varphi$ be a positive definite radial function and $A \in\mathbb{R}^{d\times d}$ be a
chosen invertible matrix. Then, the anisotropic RBF interpolant $S_{A}$
is defined by:
\begin{equation}
S_{A}(\mathbf{x})= \sum^{N}_{j=1}c_{j}\varphi_{A}(\|
\mathbf{x}-\mathbf{x}_{j}\|), \quad
\mathbf{x} \in \Omega,
\end{equation}
for $c_i$ such that the interpolation conditions $S_A(\mathbf{x}_i)=y_i$, $i=1,\dots,N$, are satisfied.

The well-posedness of the interpolation problem is guaranteed (for positive definite kernels) as a direct consequence of the invertibility of the scaling matrix $A$. We refer to~\cite{Ani10} for the error analysis of anisotropic RBF
interpolation.

\section{Sparse kernel-based interpolation}\label{sec3}
We describe the basic \emph{sparse kernel-based interpolation} (SKI) method that will be used in each step of the multilevel algorithm. The SKI method can be used also as a single step interpolation method. 

The starting point is the observation that, assuming sufficient smoothness of the interpolation data, the number of points required to provide a given accuracy can be dramatically reduced when basis functions with carefully constructed direction-wise anisotropic scaling are used. This way, notwithstanding the approximation strength coming from a few directions only, there is only negligible loss of accuracy, due to the additional smoothness assumed. Hyperbolic cross products and sparse grids use this idea in the context of high dimensional numerical quadrature, approximation and in numerical solution of partial differential equations. 

In~\cite{SchreiberAnja2000}, hyperbolic crosses of tensor products of one-dimensional RBFs have been considered without numerical assessment of the resulting method. The use of non-tensor product $d$-dimensional RBFs in the hyperbolic-cross/sparse-grid setting is not straightforward, as such approximation spaces are characterised by basis functions with differently anistropic scalings in different directions. Using such scalings would \emph{not} guarantee the well-posedness of the resulting kernel-based interpolation problem.

An alternative approach can be motivated by the, so-called, $d$-boolean interpolation \cite{DelvosFJ1982} or, the related combination technique \cite{GSZ92,GarckeandGriebel2001}, for piecewise polynomial interpolation on sparse grids. A key point in this approach is that the piecewise polynomial interpolant is equivalent to a linear combination of interpolants constructed on a carefully selected collection of subsets of the (anisotropic) basis functions. Each such subset consists of translations of the
 \emph{same} (anisotropic) basis function.

To construct the sparse kernel-based interpolant, we solve a number of anisotropic radial basis function interpolation problems on appropriately selected sub-grids and we linearly combine the resulting partial interpolants to obtain the sparse kernel-based interpolant.

More specifically, let $ \Omega$ $ \mathrel{\mathop:}= $ $[ 0,1 ]^{d}$, and let
$u:\Omega \rightarrow
 \mathbb{R}$. The extension to general axiparallel domains is straightforward; we refrain from doing so here in the interest of simplicity of the presentation only. Comments on possible extensions to more general domains will be given in Section \ref{sec7}. 
 
For a multi-index $\mathbf{l} = (l_{1},\ldots,l_{d})\in
\mathbb{N}^{d}$,  we define the family of directionally uniform grids
$\{\mathbb{X}_{\mathbf{l}}:\mathbf{l}\in
\mathbb{N}^{d}\}$, 
in $\Omega$, with meshsize
$h_{\mathbf{l}}=2^{-\mathbf{l}}:=(2^{-l_{1}},\ldots,2^{-l_{d}})$. That is, 
 $\mathbb{X}_{\mathbf{l}}$ consists of the points
 $ \mathbf{x}_{\mathbf{l},\mathbf{i}}:= (x_{l_{1},i_{1}},\ldots,x_{l_{d},i_{d}}) $,
 with $x_{l_{j},i_{j}}=i_{j}2^{-l_{j}}$, for $ i_{j}=0,1,\dots, 2^{l_{j}}$, $j=1,\dots, d$.
 The number of nodes
$N^{\mathbf{l}}$ in $\mathbb{X}_{\mathbf{l}}$ is given by
$$
N^{\mathbf{l}}=\prod_{i=1}^{d}(2^{l_{i}}+1).
$$
If $h_{l_{i}}$ = $2^{-n},$  for all $i=1,\cdots,d,$
 $\mathbb{X}_{\mathbf{l}}$ is the uniform \emph{full grid} of level
$n$, having
 size $N$=$(2^{n}+1)^{d}$; this will be denoted by $\mathbb{X}^{n,d}$.
 
\begin{figure}
\begin{center}
\begin{tabular}{m{1.8cm}m{0.1cm}m{1.6cm}m{0.1cm}m{1.6cm}m{0.1cm}m{1.6cm}m{0.1cm}m{1.6cm}}
   \includegraphics[width=2.2cm]{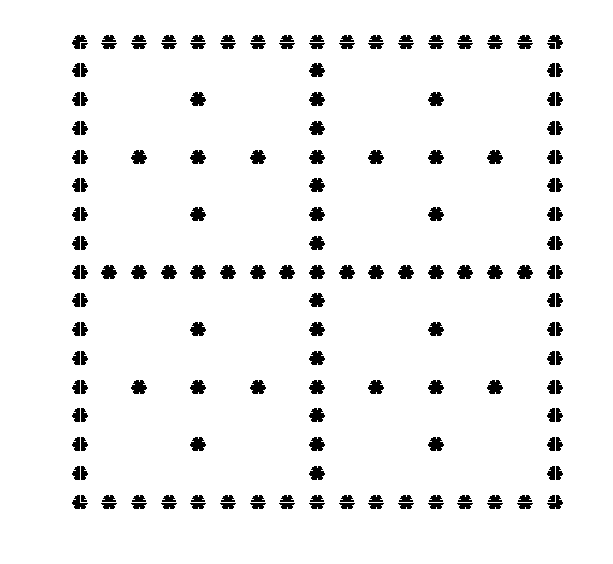} &=
      &\includegraphics[width=2cm]{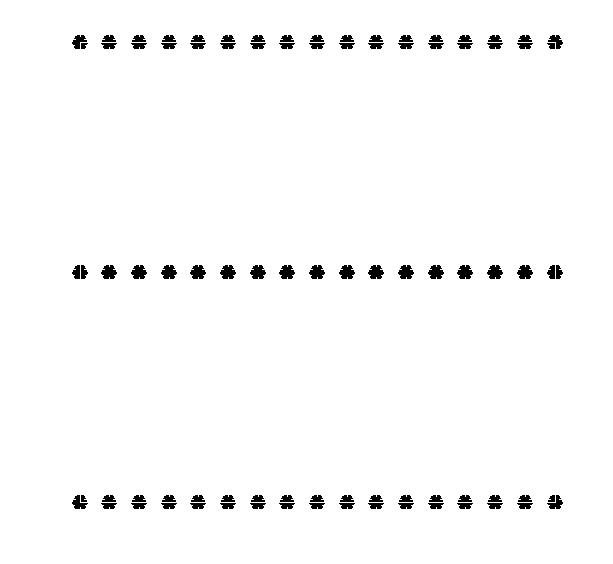} 
     &$\cup$  &\includegraphics[width=2cm]{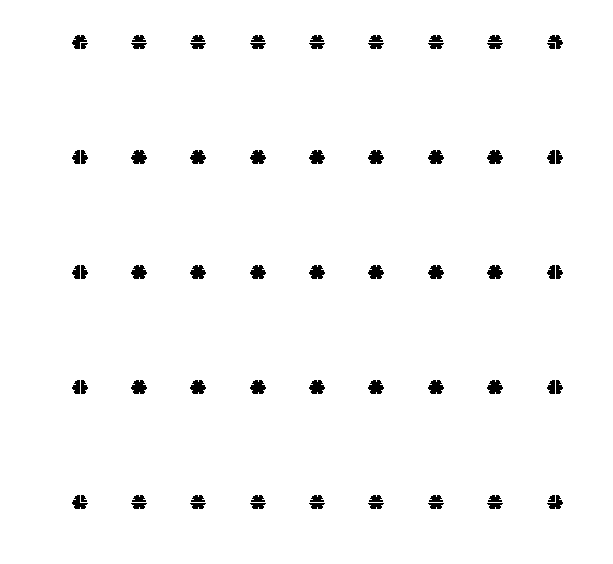} 
     & $\cup$ & \includegraphics[width=2cm]{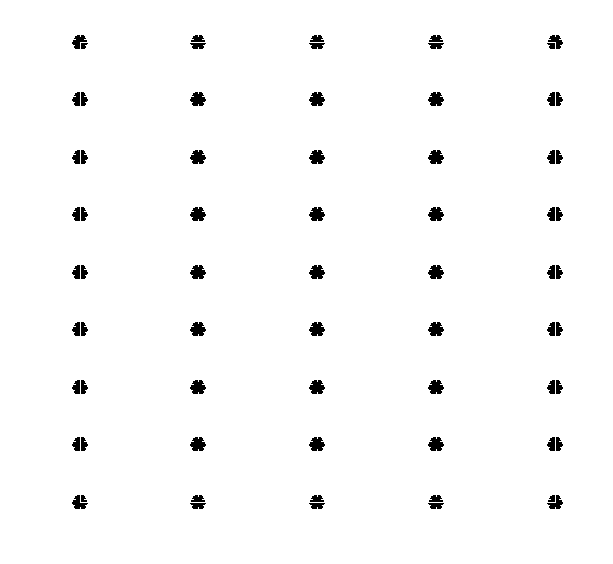} 
     & $\cup$ & \includegraphics[width=2cm]{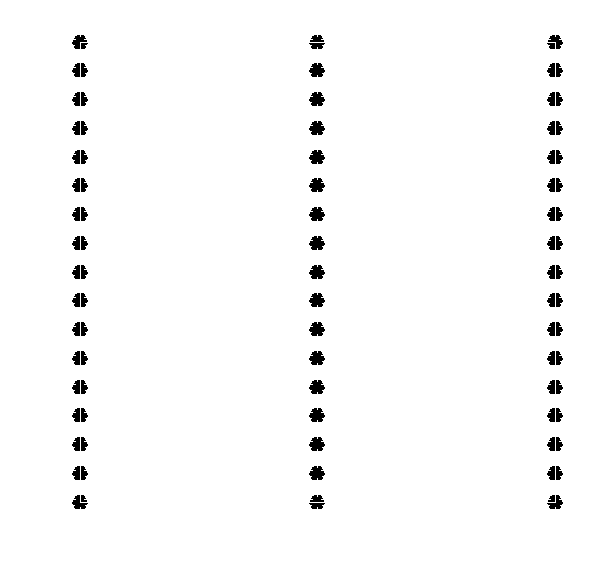}  
\end{tabular}
\end{center}
\caption{Sparse grid $\tilde{\mathbb{X}}^{4,2}$ via \eqref{sparse_grid_def}. } \label{Fig:Sparse grid def}
\end{figure}

We also consider the following subset of $\mathbb{X}^{n,d}$,
\begin{equation}\label{sparse_grid_def}
\tilde{\mathbb{X}}^{n,d}:=\bigcup_{ | \mathbf{l}|_{1} = n+(d-1)}\mathbb{X}_{\mathbf{l}},
\end{equation}
with $|\mathbf{l}|_1:=l_1+\dots+l_d$, which will be referred to as the \emph{sparse grid of level $n$ in $d$ dimensions}.  We refer to Figure \ref{Fig:Sparse grid def} for a visual representation of \eqref{sparse_grid_def} for $n=4$ and $d=2$. Notice that there is some redundancy in this definition as some grid points are included in more than one sub-grid.

We want to evaluate the interpolant  at the constituent sub-grids
$\mathbb{X}_{\mathbf{l}}$. As the constituent grids admit different density in each coordinate direction, we shall make use the anisotropic RBFs from Section \ref{sec2}. To this end, for each multi-index $\mathbf{l}=(l_{1},\ldots,l_{d})$, we define the transformation matrix $A_{\mathbf{l}}\in \mathbb{R}^{d \times d}$ by 
\[
A_{\mathbf{l}}:={\rm diag}(2^{l_{1}},\ldots,2^{l_{d}}).
\]

The anisotropic RBF interpolant $S_{A_{\mathbf{l}}}$ of $u$ at the points of $\mathbb{X}_{\mathbf{l}}$ is then defined by
\begin{equation}
S_{A_{\mathbf{l}}}(\mathbf{x}):=
\sum^{N^{\mathbf{l}}}_{j=1}c_{j}\varphi(\|
A_{\mathbf{l}}(\mathbf{x}-\mathbf{x}_{j})\|),
\end{equation}
for $\mathbf{x}\in\Omega$, where $c_j\in\mathbb{R}$ are chosen so that the interpolation conditions
\[
 S_{A_{\mathbf{l}}}|_{\mathbb{X}_{\mathbf{l}}}=u|_{\mathbb{X}_{\mathbf{l}}},
\]
are satisfied.

\begin{figure}
\begin{center}
\begin{tabular}{m{1.8cm}m{0.1cm}m{1.6cm}m{0.1cm}m{1.6cm}m{0.1cm}m{1.6cm}m{0.1cm}m{1.6cm}}
   \includegraphics[width=2.2cm]{gridS4.png} &=
      &\includegraphics[width=2cm]{gridX41.png} 
     & $\oplus$  &\includegraphics[width=2cm]{gridX32.png} 
     & $\oplus$ & \includegraphics[width=2cm]{gridX23.png} 
     & $\oplus$ & \includegraphics[width=2cm]{gridX14.png}  \\
     & $\ominus$  &\includegraphics[width=2cm]{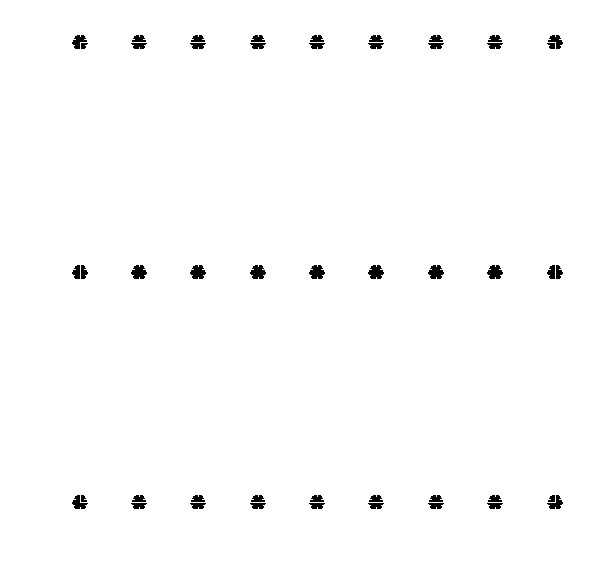} 
     & $\ominus$ & \includegraphics[width=2cm]{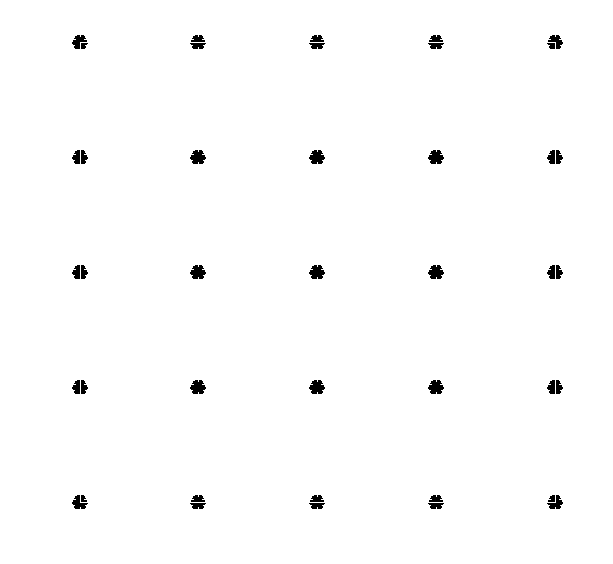} 
     & $\ominus$ & \includegraphics[width=2cm]{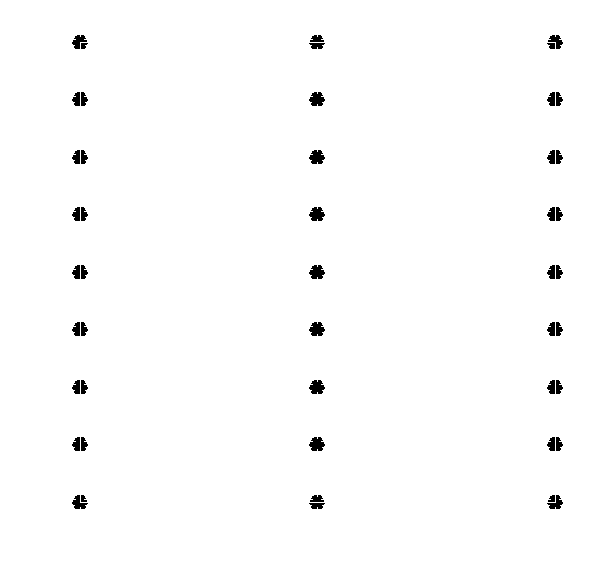}  &
\end{tabular}
\end{center}
\caption{The construction of the sparse kernel-based $\tilde{S}_{4}$ interpolant on $\tilde{\mathbb{X}}^{4,2}$.} \label{Fig:Sparse grid def2}
\end{figure}

To construct the \emph{sparse kernel-based interpolant} (SKI, for short) $S^{c}_{n} $ on the sparse grid $\tilde{\mathbb{X}}^{n,d}$, the sub-grid interpolants $S_{A_{\mathbf{l}}}$ are linearly combined using the formula
\begin{equation}\label{Eq:SPGDinterpcombination}
\tilde{S}_{n}(\mathbf{x}) =\sum^{d-1}_{q=0}(-1)^{q} \left(
                                              \begin{array}{c}
                                                d-1 \\
                                                q \\
                                              \end{array}
                                            \right)
 \sum_{| \mathbf{l}|_{1}=n+(d-1)-q} S_{A_{\mathbf{l}}}(\mathbf{x}).
\end{equation}
The combination formula~\eqref{Eq:SPGDinterpcombination} has been used in the context of $d$-boolean Lagrange polynomial interpolation~\cite{DelvosFJ1982}, and in the combination technique for the numerical solution of elliptic partial differential equations using the finite element method on sparse grids ~\cite{GSZ92,GarckeandGriebel2001}. For $d=2$, \eqref{Eq:SPGDinterpcombination} becomes
\begin{equation}\label{Eq:SPGDinterpcombination2d}
\tilde{S}_{n}(\mathbf{x}) =
 \sum_{| \mathbf{l}|_{1}=n+1} S_{A_{\mathbf{l}}}(\mathbf{x})-  \sum_{| \mathbf{l}|_{1}=n} S_{A_{\mathbf{l}}}(\mathbf{x}),
\end{equation}
that is, the first term on the right-hand side of \eqref{Eq:SPGDinterpcombination2d} gives the sum of interpolants on the sub-grids of level $n+1$, while the second term on the right-hand side  of \eqref{Eq:SPGDinterpcombination2d} subtracts the redundant points visited more than once. We refer to Figure \ref{Fig:Sparse grid def2} for an illustration when $d=2$ and $n=4$.

Hence, SKI uses a dimension-wise multilevel decomposition
of interpolation data sites in conjunction with the application of
kernel-based interpolants $S_{A_{\mathbf{l}}}(\cdot)$ with different
scaling in each direction. The sparse kernel-based interpolant $\tilde{S}_{n}$ can be implemented in a quite straightforward fashion by utilising existing, fast, RBF interpolation algorithms: the only modification needed is the introduction of a scaling for each sub-grid problem. We note that each interpolation problem can be solved completely independently, rendering the resulting SKI method ideally suited for implementation in parallel computers.

Moreover, we observe the size of each sub-grid problem is $O(2^n)$, where $n$ is the number of levels, i.e., it is \emph{independent} of the dimension $d$. One has to solve $O(n^{d-1})$ such sub-grid problems to obtain the sparse kernel-based interpolant. Observe that $O(2^n)$ is the number of points of the corresponding full grid $\mathbb{X}^n$ in \emph{each} space direction. Hence, if the rate of convergence for the SKI algorithm is comparable with the one of the standard RBF algorithm for a sufficiently large class of underlying functions $u$, the benefits in complexity are potentially significant, especially for $d\ge 2$. Indeed, when the underlying function $u$ is assumed to admit sufficient regularity (e.g., if the regularity of $u$ is characterised by anisotropic Sobolev spaces) \cite{DelvosFJ1982,GSZ92,GarckeandGriebel2001} and \cite{SchreiberAnja2000} the rate of convergence on sparse grids for tensor-product of piece-wise polynomials and of one-dimensional RBFs, respectively, is shown to be optimal (modulo a logarithmic factor). Although, there is no general proof at this point that this is also the case for the SKI method presented here,  numerical experiments presented below show that for the multilevel version of the SKI algorithm, described in the next section, good convergence results are also observed. We finally remark that for the case of Gaussian interpolation, due its tensor-product nature, the theoretical results of \cite{SchreiberAnja2000} on tensor-product one-dimensional sparse RBF interpolation should be applicable to the case of SKI also.

\section{Multilevel Sparse Kernel-Based Interpolation}\label{sec4}
 Multilevel methods for RBFs~\cite{Floater&Iske1996,Iske2001,Wendland1998,Fasshauer&Jerome1999,
Fasshauer99b,Narcowich&Schaback&Ward1999,Hales&Levesley2002,
IskeANDLeveslely2005,Ohtake2005} combine the advantages of stationary and non-stationary interpolation, aiming to accelerate convergence and to improve the numerical stability of the interpolation procedure.  The basic idea of multilevel methods in this context
is to interpolate the data at the coarsest level, and then update by interpolating the residuals on progressively finer data sets using appropriately scaled basis functions.

The setting of SKI is naturally suited to be used within a multilevel interpolation algorithm. Indeed, the sparse grids from lower to higher level are nested, i.e., $\tilde{\mathbb{X}}^{n,d}\subset \tilde{\mathbb{X}}^{n+1,d}$ for $n\in\mathbb{N}$; we refer to Figure~\ref{Fig:First
six sparse grids in 2D} for an illustration when $d=2$. Moreover,  each sub-grid interpolant  uses appropriately scaled anisotropic basis function with the scaling being proportional to density of the corresponding constituent sub-grid. Finally, due to the geometrical progression in the problem size from one sparse grid to the next, a multilevel algorithm would not affect adversely the attractive complexity properties of SKI.

The \emph{multilevel SKI} (MLSKI, for short) algorithm is initialised by computing the SKI $\tilde{S}_{n_0^{}}$ at the coarsest designated sparse grid $\tilde{\mathbb{X}}^{n_0^{},d}$ and set $\Delta_0 :=\tilde{S}_{n_0^{}}$. Then, for $k=1,\dots n$,  we compute $\Delta_k$ to be the sparse grid interpolant of the residual $u-\sum_{j=0}^{k-1}\Delta_j$ on $\tilde{\mathbb{X}}^{k,d}$. The resulting multilevel sparse kernel based interpolant is then given by 
\[
\tilde{S}^{\rm ML}_n : = \sum_{j=0}^{n}\Delta_j.
\]
\begin{figure}
\begin{center}
\includegraphics[width=8.5cm]{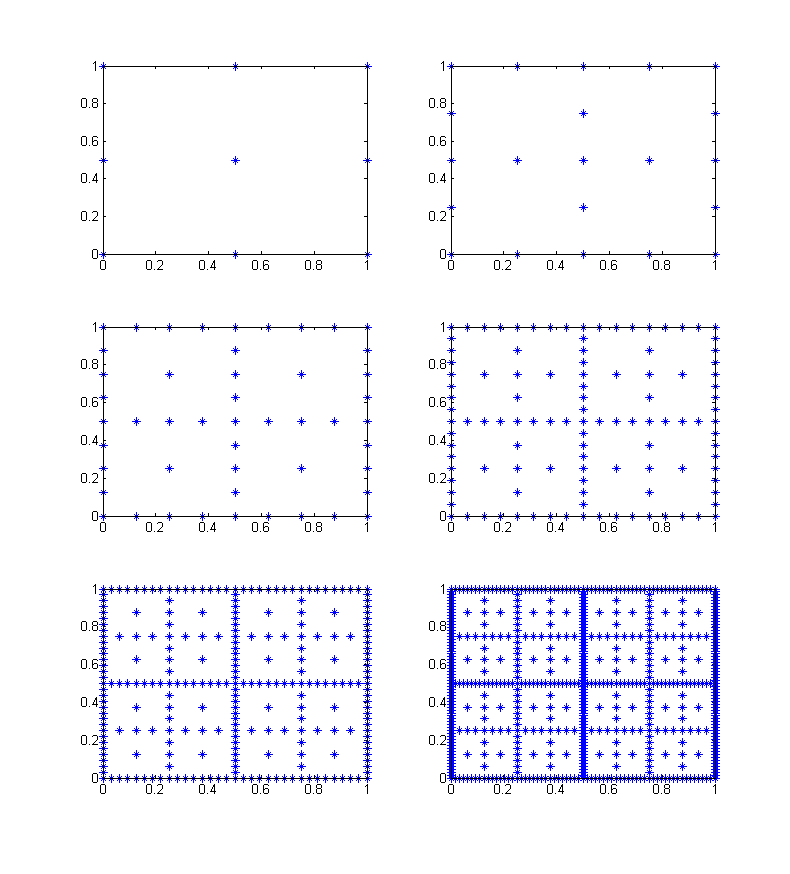}
\end{center}
\caption{Nested-ness of sparse grids $\tilde{\mathbb{X}}^{n,2}$, $n=1,2,\dots, 6$.}
 \label{Fig:First six sparse grids in 2D}
\end{figure}

\section{Stability of SKI and MLSKI}\label{sec5}
The numerical stability of RBF interpolation is a challenging issue, due to the ill-conditioning of the respective interpolation matrices in standard bases; we refer to~\cite{Wend05, Fasshauerbook2007} and the references therein for a discussion. It is, therefore, evident that
the numerical stability of SKI and MLSKI algorithms can be
assessed by considering the maximum condition number of the constituent sub-grid
interpolation problems.

The main factors  affecting the conditioning of an RBF interpolation problem are the problem size, the separation distance $q_X$ of the respective data set $X$, given by
\[
q_{X}:= \frac{1}{2}\min_{i,j}\|{\mathbf x}_i -{\mathbf x}_j \|_2, \quad {\mathbf x}_i, {\mathbf x}_j\in X,
\]
 and the shape parameters of the RBF (which adjusts the ``strength'' or the support, depending on the type of RBF). The choice of shape parameter and the monitoring of the separation distance require more attention in practice than the problem size, which is of secondary importance~\cite{Wend05, Fasshauerbook2007}.
 Recipes for choosing RBF shape parameters in various contexts have been presented, e.g., in~\cite{Hardy71,Franke1982b,CarlsonandFoley1991, CarlsonandFoley1992, KansaandCarlson1992,Carlson-and-Natarajan1994}, while the
dependence of the stability on the shape has been addressed
in~\cite{RippaShmuel999,FB04,FP07,fasshauer_zhang,FasshauerandMccourt2011,BEN}, where stable algorithms for RBF interpolation with small
shape parameters have been proposed. 

SKI and MLSKI naturally include scaling adjustment at each level. So, the condition number of each sub-grid interpolation matrix appears to be (mildly) affected by the sub-grid size, as we shall see in Section \ref{sec6}.  This is, indeed, a very attractive feature of the SKI and MLSKI algorithms as the growth of degrees of freedom on each sub-grid grows \emph{independently} of the problem dimension $d$. 

Observing that, due to the anisotropy of the scaling matrices $A_{\mathbf l}$, the separation distance on each pulled-back (multiplied by $A_{\mathbf l}^{-1}$) sub-grid is constant with respect to the level $n$. This suggests that we can use standard heuristics for the choice of shape-parameters for multilevel RBF interpolation (see \cite{Floater&Iske1996,Iske2001,Wendland1998,Fasshauer&Jerome1999,
Fasshauer99b,Narcowich&Schaback&Ward1999,Hales&Levesley2002,
IskeANDLeveslely2005,Ohtake2005}, or \cite{Wend05} for a review). In particular,  we set
\begin{equation}\label{shape_choice}
c=\frac{q_{\mathbb{X}^{n+1,d}_{}}}{Kq_{\mathbb{X}^{n,d}_{}}},
\end{equation}
for some constant $K>0$. Extensive numerical experiments presented in \cite{Fazli_thesis} suggest that the choice $K=1$, gives very small condition numbers at the cost of (mildly) lower accuracy of the interpolation problem. The choice
$K=3$, produces larger but mostly safe condition numbers (in this work we shall refer to a condition number $\kappa$ as being \emph{safe} when $\kappa\lesssim 10^{10}$), with the resulting computations admitting good convergence properties. For $K=3$, we have $0.2\le c\le 0.8$. For $K> 3$, the ill-conditioning gradually increases with $K$, resulting in unstable computations.

\section{Numerical Experiments}\label{sec6}
We present a collection of numerical experiments  for $d=2,3,4$, where the implementation of SKI and MLSKI algorithms is assessed and compared against both the standard RBF interpolation on \emph{uniform full grids} and its standard multilevel version (MLRBF); see, e.g, \cite{Wend05} for a review.

The algorithm implementation has been performed in MATLAB$^\copyright$ on a quad-core 2.67GHz Intel Xeon X5550 CPU with 12GB of RAM, \emph{without} taking advantage of the possibility of parallel implementation of SKI or MLSKI. The SKI/MLSKI algorithm is able to solve $d$-variate interpolation problems on sparse grids up to $114,690$ centers for $d = 3$, and up to $331,780$ centers for $d = 4$, respectively. The same basic interpolation solver for classical RBF on the same machine could only solve problems of size nearly $15,000$ regardless of the dimension $d$, due to the size of the resulting linear system.

We stress that no attempt has been made to use fast algorithms, e.g., the fast Gauss transform of Greengard and Strain \cite{greengard} for the RBF interpolation problems for either the standard RBF, MLRBF or the SKI, MLSKI algorithms. The use of fast or stable algorithms for RBF interpolation could be also used within the sub-grid interpolation problems of SKI/MLSKI. Indeed, any modern fast evaluation or numerical stabilization technique in RBF interpolation of a single problem can be incorporated into the SKI/MLSKI framework, possibly substantially improving the results presented here.

\subsection{Experiment 1}
\begin{table}
\begin{center}
\addtolength{\tabcolsep}{-1pt}
\begin{tabular}{||c|c|c|c|c||}
  \hline
SGnode  &Max-error &  RMS-error &  Cond. no  & Time\\
\hline \hline
  9          & 6.2215e-1 & 1.8363e-1   & 2.6912e+3 &  $<$1\\
  21         & 3.3237e-1 & 7.6547e-2   & 2.5325e+4 & $<$1\\
  49        & 1.1130e-1 & 3.8660e-2   & 2.8184e+5 & $<$1\\
  113      & 4.0379e-2 & 1.0835e-2   & 2.6522e+6 & $<$1\\
  257      & 1.2649e-2 & 2.5117e-3   & 2.9516e+7 &  $<$1\\
  577     & 2.4678e-3 & 4.0273e-4  & 1.7591e+8 &  $<$1\\
  1281   & 2.2043e-4 & 2.1030e-5  & 1.0484e+9 &  1\\
  2817    & 3.5287e-5 & 2.5391e-6   & 2.3229e+9 &  5\\
  6145    & 6.2139e-6 & 3.2696e-7   & 5.1468e+9 &  33\\
  13313   & 1.1784e-6 & 4.2920e-8   & 6.5016e+9 &  281\\
  28673   & 2.1204e-7 & 5.6557e-9   & 8.2129e+9 &  2397\\
  61441  & 4.1321e-8 & 7.6854e-10  & 8.7056e+9 &  21586\\[0.5ex]
 \hline
\end{tabular}
\end{center}
\caption{MLSKI for $u_{F2D}$, $d=2$, using Gaussians with
$c=0.45$. Error evaluated at
$25,600$  Halton points.}
\label{Tab:MLSR2Dlevel12}
\end{table}

The MLSKI algorithm with Gaussian basis functions is applied to the (standard) benchmark Franke's function $u_{F2D}$ for $d=2$, given by
\begin{equation}
\begin{aligned}
 u_{F2D}(x_{1},x_{2})  := & \frac{3}{4}e^{(-(9x_{1}-2)^{2}-(9x_{2}-2)^{2}}
     +\frac{3}{4}e^{-((9x_{1}+1)^{2})/49-((9x_{2}+1)^{2})/10}\\
      & +\frac{1}{2}e^{-((9x_{1}-7)^{2})/4-(9x_{2}-3)^{2}}
      -\frac{1}{5}e^{-((9x_{1}-4)^{2})/4-(9x_{2}-7)^{2}}.
\end{aligned}
\end{equation}
In Table \ref{Tab:MLSR2Dlevel12}, ``SGnode" stands for the number
of sparse grid centers used, ``Max-error'' and ``RMS-error' are the maximum norm and root-mean-square ($L^2$-norm) errors, respectively, evaluated at $25,600$ Halton points \cite{kuipers}, ``Cond.~no" stands for the largest $2$-norm condition number of the sub-grid interpolation matrices. Finally, ``Time'' is the time in seconds required to solve the interpolation problem in the above computer. We note that the {\tt tic;} and {\tt toc;} commands of MATLAB$^\copyright$ have been used to produce the run times, which is accurate for longer run times.

In Table \ref{Tab:MLSKI3D-Alice}, the MLSKI algorithm with Gaussian basis functions is applied to a three-dimensional version of Franke's function $u_{F3D}$, given by

\begin{equation}
\begin{aligned}
 u_{F3D}(x_{1},x_{2},x_{3}) = & \frac{3}{4}e^{(-(9x_{1}-2)^{2}-(9x_{2}-2)^{2} -(9x_{3}-2)^{2})/4}\\
    &  +\frac{3}{4}e^{-((9x_{1}+1)^{2})/49-((9x_{2}+1)^{2})/10-((9x_{3}+1)^{2})/29}\\
      & +\frac{1}{2}e^{-((9x_{1}-7)^{2})/4-(9x_{2}-3)^{2}-((9x_{3}-5)^{2})/2}\\
       & -\frac{1}{5}e^{-((9x_{1}-4)^{2})/4-(9x_{2}-7)^{2}-((9x_{3}-5)^{2})}.
\end{aligned}
\end{equation}

\begin{table}
\begin{center}
\addtolength{\tabcolsep}{-1pt}
\begin{tabular}{||c|c|c|c|c||}
  \hline
SGnode  & Max-error &  RMS-error &  Cond. no  & Time\\
\hline \hline
  27        & 6.8808e-1 & 1.0179e-1     & 1.4863e+4   &   $<$1 \\
  81         & 5.5853e-1   & 7.7339e-2  & 2.5376e+5   &   $<$1 \\
  225       & 2.4324e-1   & 3.8389e-2  & 6.2956e+4   &   $<$1 \\
  593       & 1.5884e-1   & 2.1676e-2  & 1.0617e+7   &  8 \\
  1505     & 6.2918e-2   & 6.7591e-3  & 6.9125e+5   &  9 \\
  3713     & 1.3403e-2   & 1.7755e-3  & 4.2986e+8   &  12 \\
  8961     & 2.2041e-3   & 2.2448e-4  & 6.2618e+6   &  42 \\
  21249    & 3.3081e-4  & 2.9755e-5  & 9.2544e+9   &  328 \\
  49665    & 8.9456e-5  & 4.5151e-6 & 2.1038e+7    &  2900 \\
  114689  & 1.5829e-5  & 5.7471e-7 & 1.0594e+11   &  25160 \\[0.5ex]
  \hline
\end{tabular}
\end{center}
\caption{MLSKI for $u_{F3D}$, using Gaussians with
$c$ as in (\ref{shape_choice}) for $K=3$. Error evaluated at
$125,000$  Halton points.}
\label{Tab:MLSKI3D-Alice}
\end{table}

Next, we consider the interpolation problem of the simple function $u_{quad}:[0,1]^4\to\mathbb{R}$, with
\[
u_{quad}(x_1,x_2,x_3,x_4):= 4^4 \prod_{i=1}^4x_i(1-x_i).
\]
 We employ the SKI algorithm for $d=4$; the results are given in Table \ref{Tab: SKI-4D_alice}. 
 
\begin{table}
\begin{center}
\addtolength{\tabcolsep}{-1pt}
\begin{tabular}{||c|c|c|c|c||}
  \hline
SGnode  & Max-error &  RMS-error &  Cond. no  & Time\\
\hline \hline
  81          & 7.9105e-2 & 4.4589e-2  & 3.6544e+5  &  $<$1\\
  297        & 2.4067e-2 & 1.0677e-2  & 8.6224e+6  &  $<$1\\
  945       & 1.9844e-2 & 6.3598e-3  & 1.0568e+6  &  $<$1 \\
  2769      & 5.6653e-3 & 1.2672e-3  & 3.6076e+8  &  2\\
  7681      & 4.7096e-3 & 8.2613e-4  & 1.4065e+7  &  14\\
  20481     & 1.3155e-3 & 1.5425e-4  & 1.4848e+10 &  164\\
  52993    & 1.1548e-3 & 1.0690e-4  & 1.2741e+8  &  1766\\
  133889   & 3.2099e-4 & 1.9243e-5  & 6.0115e+11 &  16442\\
  331777  & 2.8385e-4 & 1.3934e-5  & 1.1542e+9  &  169643\\[0.5ex]
 \hline
\end{tabular}
\end{center}
\caption{SKI for $u_{quad}$, using Gaussians with
$c$ as in (\ref{shape_choice}) for $K=3$. Error evaluated at
  $194,481$
  Halton points.} \label{Tab: SKI-4D_alice}
\end{table}

\subsection{Experiment 2} We now turn our attention to the conditioning of the sub-grid interpolation problems. In the previous experiment a, mostly safe, maximum condition number for each level's sub-grid interpolation problems is observed, using $K=3$ in \eqref{shape_choice}. Here, we investigate the choice $K=1$ which gives $O(1)$ condition numbers at the cost of a slight reduction in the quality of the approximation in the MLSKI algorithm. 

The results from the MLSKI algorithm using Gaussians with shape parameter $c$ chosen from \eqref{shape_choice} with $K=1$ are given in Table  \ref{MLSKI_FRANKE3D_K1}. The ``Time''-column is ommitted for brevity as it is almost identical to the one of Table \ref{Tab:MLSKI3D-Alice}. We observe that the error is roughly 5 times worse, but the condition numbers of the resulting sub-grid problems are nearly optimal. 

\begin{table}
\begin{center}
\addtolength{\tabcolsep}{-5pt}
\begin{tabular}{||c|c|c|c||}
  \hline
SGnode  & Max-error &  RMS-error &  Cond. no  \\
\hline\hline
  27     & 7.0968e-1&  1.0531e-1  & 1.8\\
  81     &  5.5864e-1 & 7.6444e-2  & 2.6\\
  225   & 3.2513e-1&  4.8817e-2 & 1.5\\
  593  & 1.3272e-1 & 1.5118e-2  & 3.0\\
  1505 & 7.8689e-2 & 7.9827e-3 & 1.5\\
  3713  & 2.1970e-2& 2.1392e-3 & 3.2 \\
  8961  & 1.0543e-2&  9.9965e-4 & 1.6\\
  21249  & 1.7569e-3&  1.7839e-4   &3.4\\

 \hline
\end{tabular}
\end{center}
 \caption{MLSKI for $u_{F3D}$, using Gaussians with
$c$ as in (\ref{shape_choice}) for $K=1$. Error evaluated at
$125,000$  Halton points.}
\label{MLSKI_FRANKE3D_K1}
\end{table}

\subsection{Experiment 3}

 We continue by comparing the SKI and MLSKI algorithms with standard RBF interpolation on full grids and with its (standard) multilevel version on full grids, henceforth denoted by MLRBF. 
 
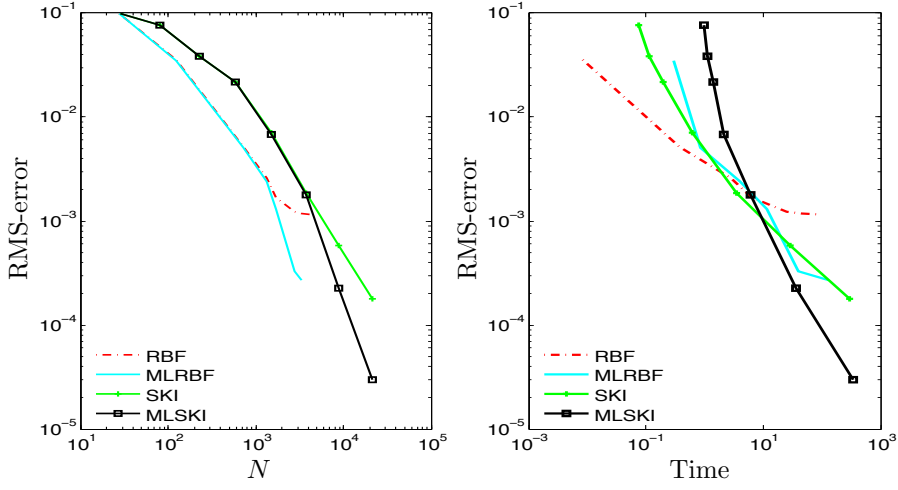
\begin{figure}
\begin{center}
\subfigure{
\includegraphics[width=5.5cm,height=6cm]{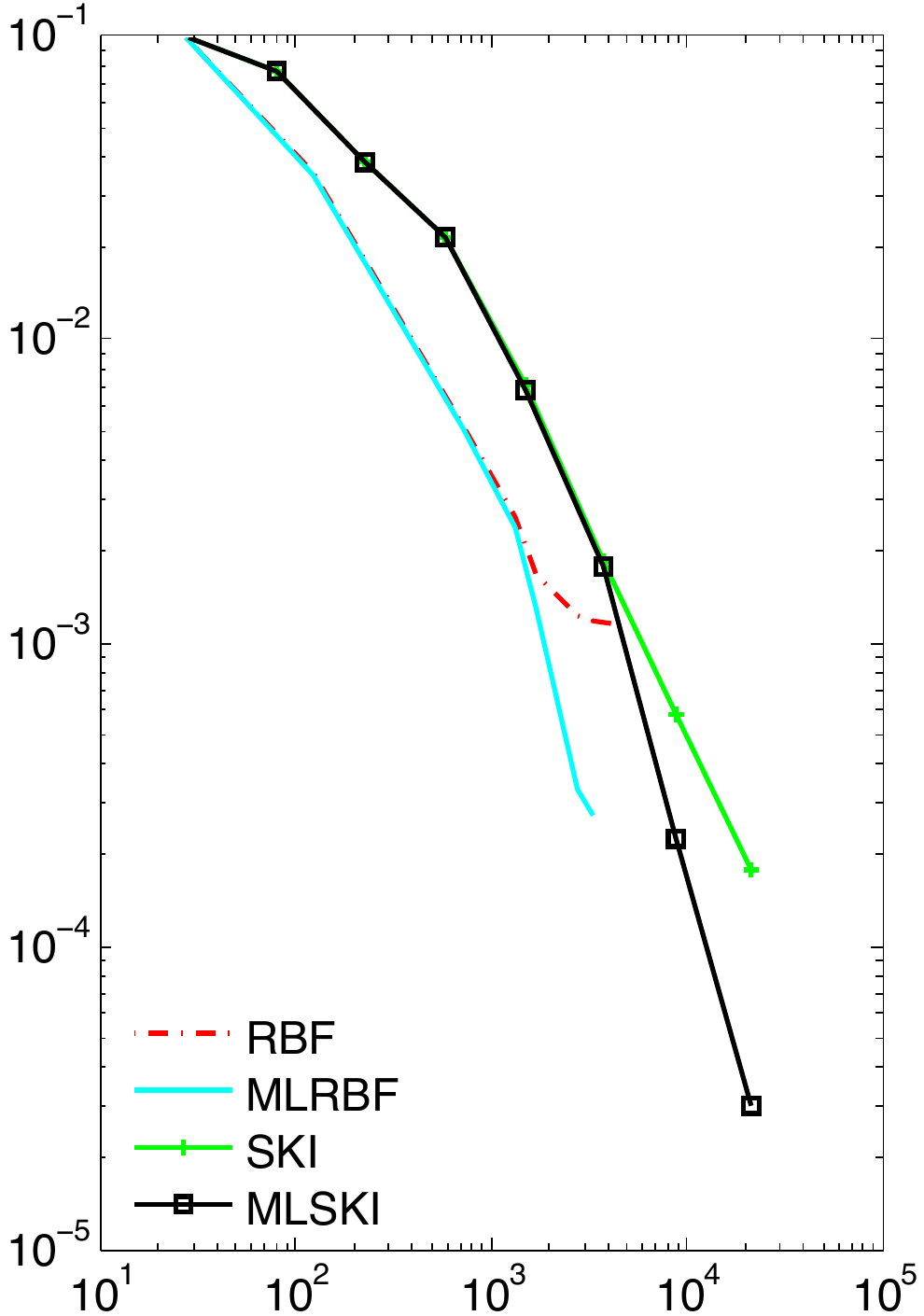}
\begin{picture}(0,0)
\put(-80,-10){$N$} \put(-170,65){\rotatebox{90}{RMS-error}}
\end{picture}
}
\subfigure{
\includegraphics[width=5.5cm,height=6cm]{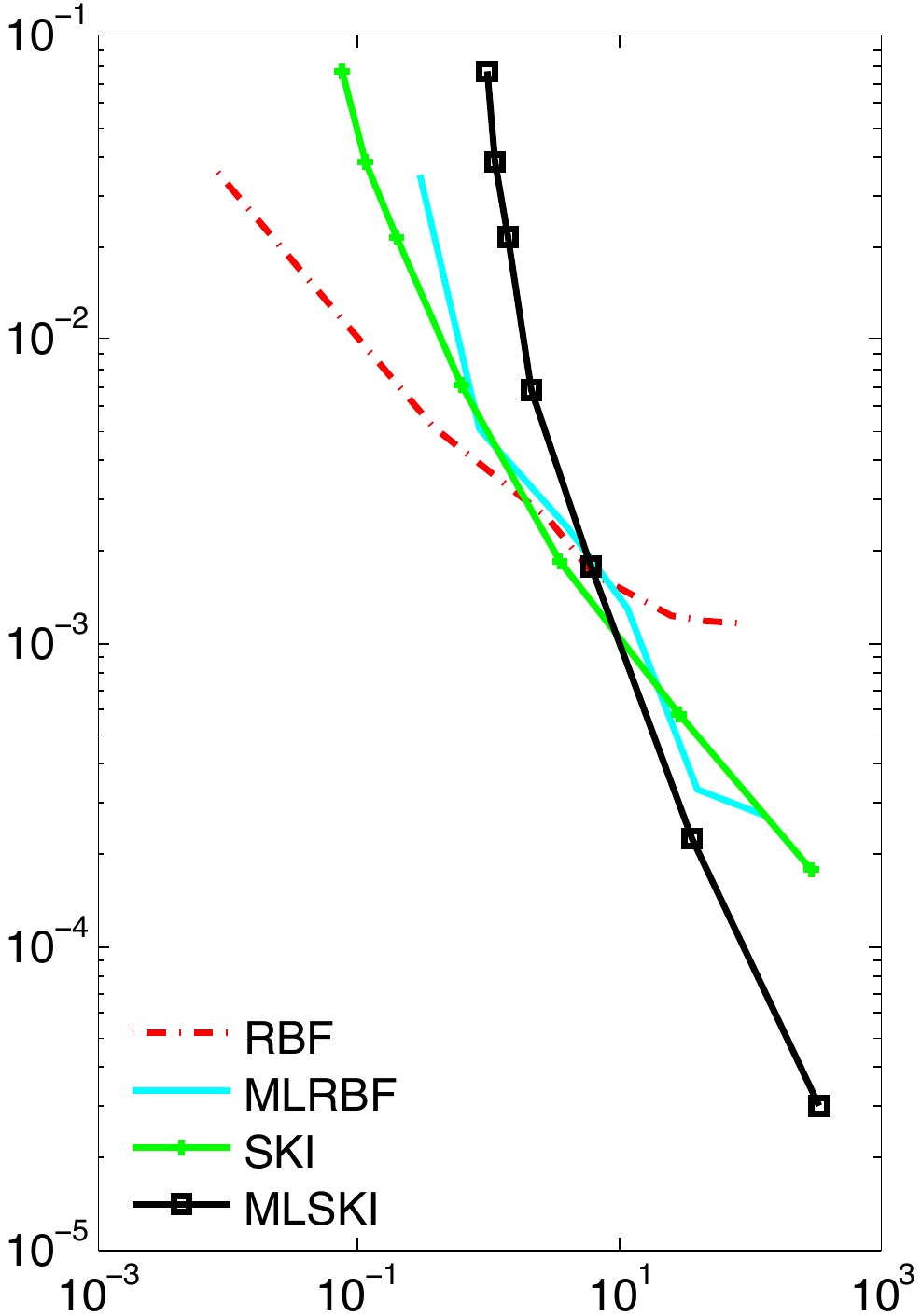}
\begin{picture}(0,0)
\put(-90,-10){Time} \put(-170,65){\rotatebox{90}{RMS-error}}
\end{picture}
}
\end{center}
\caption{Convergence of Gaussians for $u_{F3D}$. Error evaluated at
$125,000$ Halton points.} \label{F3D}
\end{figure}

In Figure \ref{F3D}, the root mean-square error of RBF, MLRBF, SKI and MLSKI, respectively,  for $u_{F3D}$, are plotted against the number of data sites $N$ and against the computation time (in seconds). For RBF and MLRBF, $N$ denotes the number of data sites on the full grid $\mathbb{X}^n$, whereas for SKI and MLSKI, $N$ refers to the number of sparse grid nodes (i.e., $N=$ SGnode). (We note in passing that using standard ``isotropic'' RBF interpolation on sparse grids does not appear to result to a convergent algorithm.) The choice of the shape parameter is given by \eqref{shape_choice} with $K=3$ for SKI/MLSKI for all experiments; the shape parameter for RBF/MLRBF is chosen so as to have comparable condition numbers in all cases. We observe that RBF/MLRBF interpolation appears to perform better than SKI/MLSKI when the error is plotted against $N$. The SKI/MLSKI algorithms are, on the other hand, able to calculate larger problems and they do so efficiently in terms of complexity, at least for the case of larger problems. This is manifested in the RMS-error versus computation time plot in Figure \ref{F3D}.

\begin{figure}
\begin{center}
\subfigure{
\includegraphics[width=5.5cm,height=6cm]{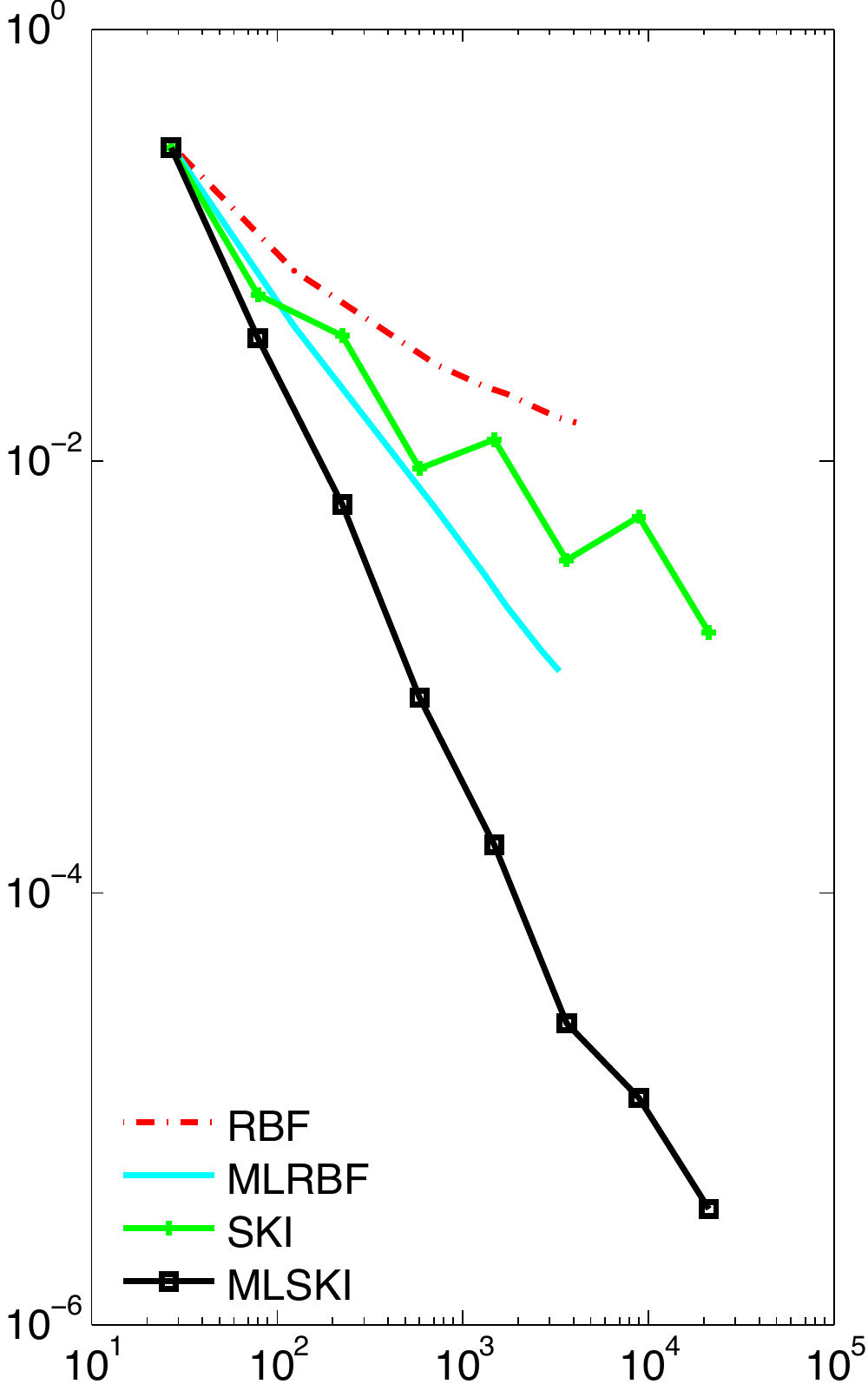}
\begin{picture}(0,0)
\put(-80,-10){$N$} \put(-170,65){\rotatebox{90}{RMS-error}}
\end{picture}
}
\subfigure{
\includegraphics[width=5.5cm,height=6cm]{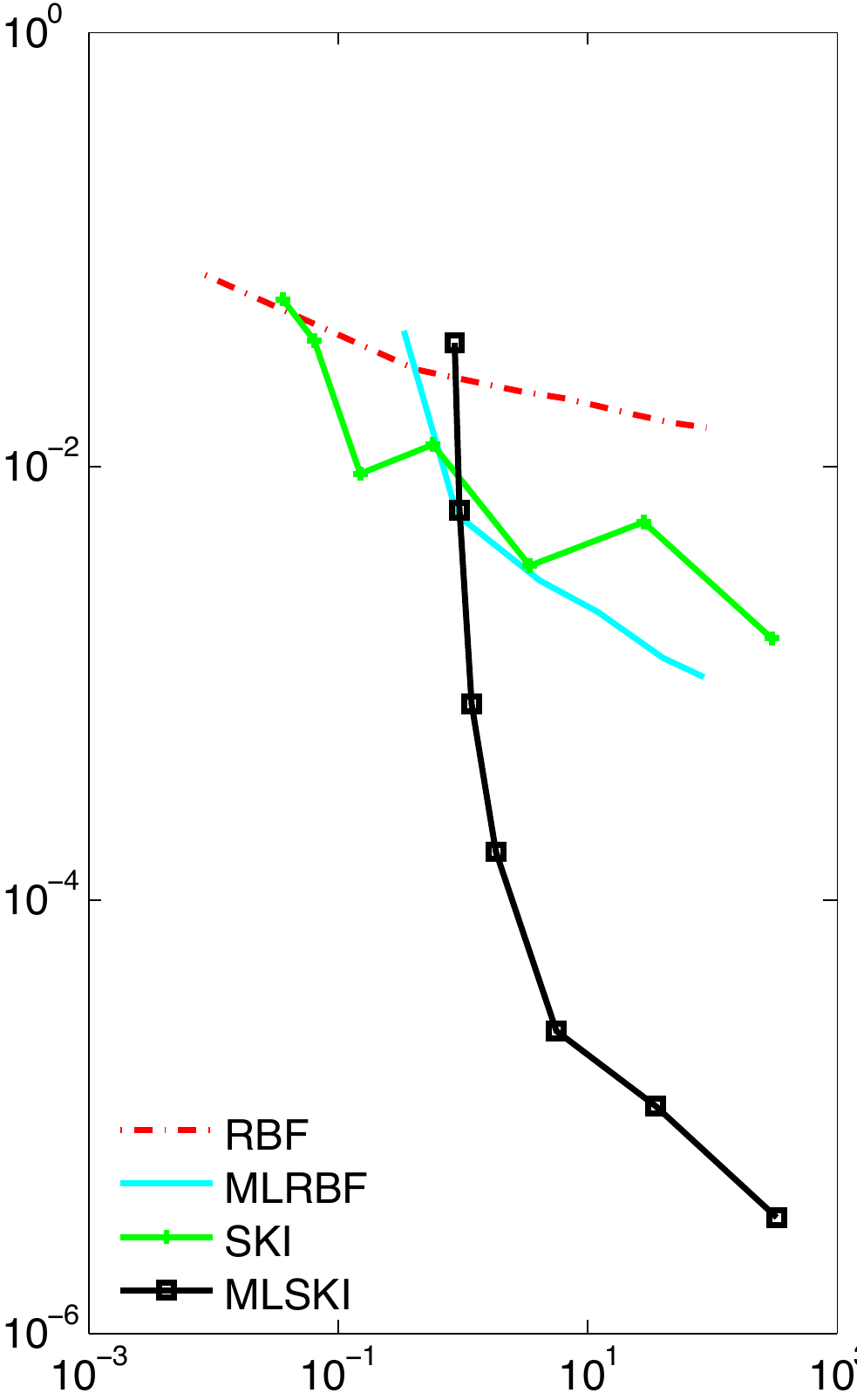}
\begin{picture}(0,0)
\put(-90,-10){Time} \put(-170,65){\rotatebox{90}{RMS-error}}
\end{picture}
}
\end{center}
\caption{Convergence of Gaussians for $u_{R3D}$. Error evaluated at
$125,000$ Halton points.} \label{R3D}
\end{figure}
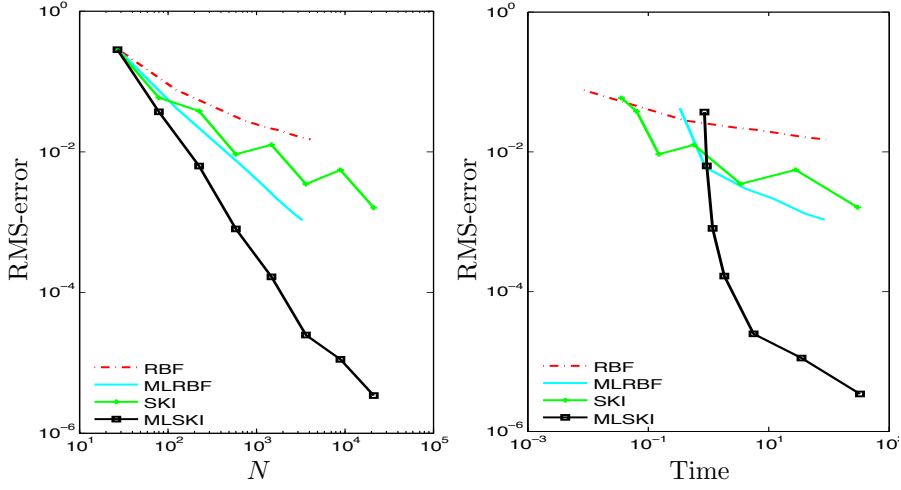

 As discussed in Section \ref{sec3}, $d$-Boolean Lagrange interpolation using polynomials requires additional smoothness of mixed derivatives in order to ensure the essentially optimal convergence rate. To test the extend to which this is also the case for the MLSKI algorithm, we consider the function $u_{R3D}:[0,1]^3\to\mathbb{R}$, with
\[
  u_{R3D}(x_{1},x_{2},x_{3}) := (r^2+r^4)\log r, \quad \text{with } \;r=(x _{1}^{2}+x _{2}^{2}+x _{3}^{2})^{1/2}.
\]
This function does \emph{not} posses smooth mixed derivatives of (arbitrarily) high order.
In Figure \ref{R3D}, the root mean-square error of RBF, MLRBF, SKI and MLSKI, respectively,  for $u_{R3D}$, are plotted against the number of data sites $N$ and against the computation time (in seconds). Somewhat surprisingly, we observe that  the MLSKI algorithm outperforms RBF/MLRBF interpolation both in terms of degrees of freedom and in terms of computational time, at least for larger problems.

Next, we turn our attention to the problem of interpolating five-dimensional data, i.e., $d=4$. This is a challenge in practice with known methods: standard RBF interpolation is limited to few tens of thousands of data points. This, in turn, results to using $O(10)$ data points in each coordinate direction, thereby, limiting the resulting approximation quality.  

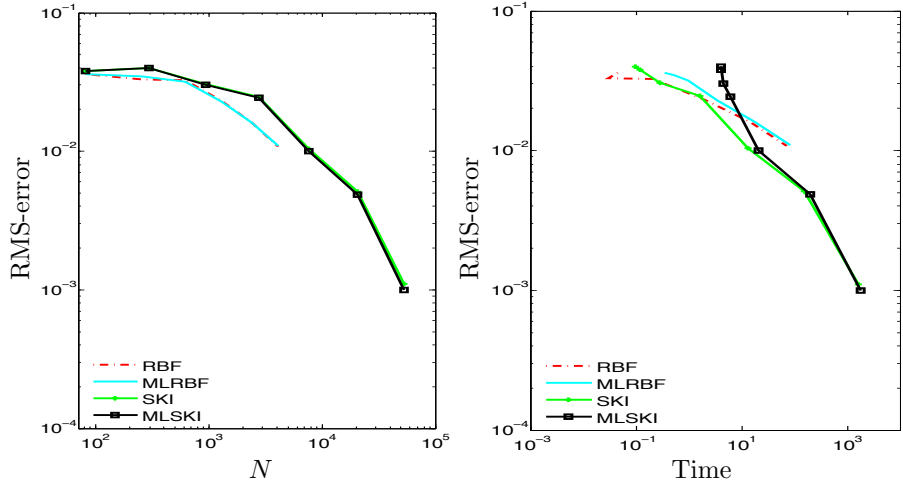
\begin{figure}
\begin{center}
\subfigure{
\includegraphics[width=5.5cm,height=6cm]{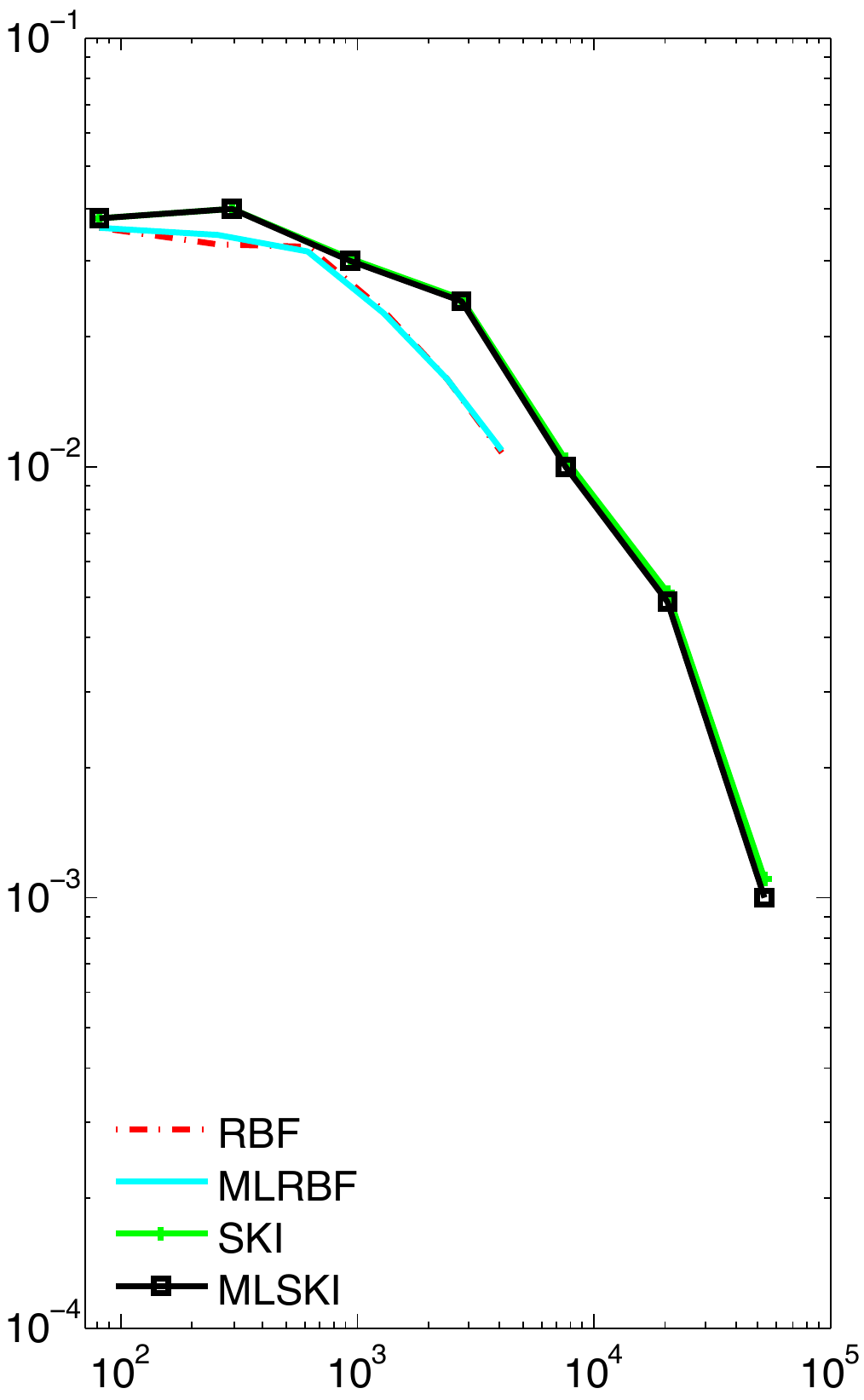}
\begin{picture}(0,0)
\put(-80,-10){$N$} \put(-170,65){\rotatebox{90}{RMS-error}}
\end{picture}
}
\subfigure{
\includegraphics[width=5.5cm,height=6cm]{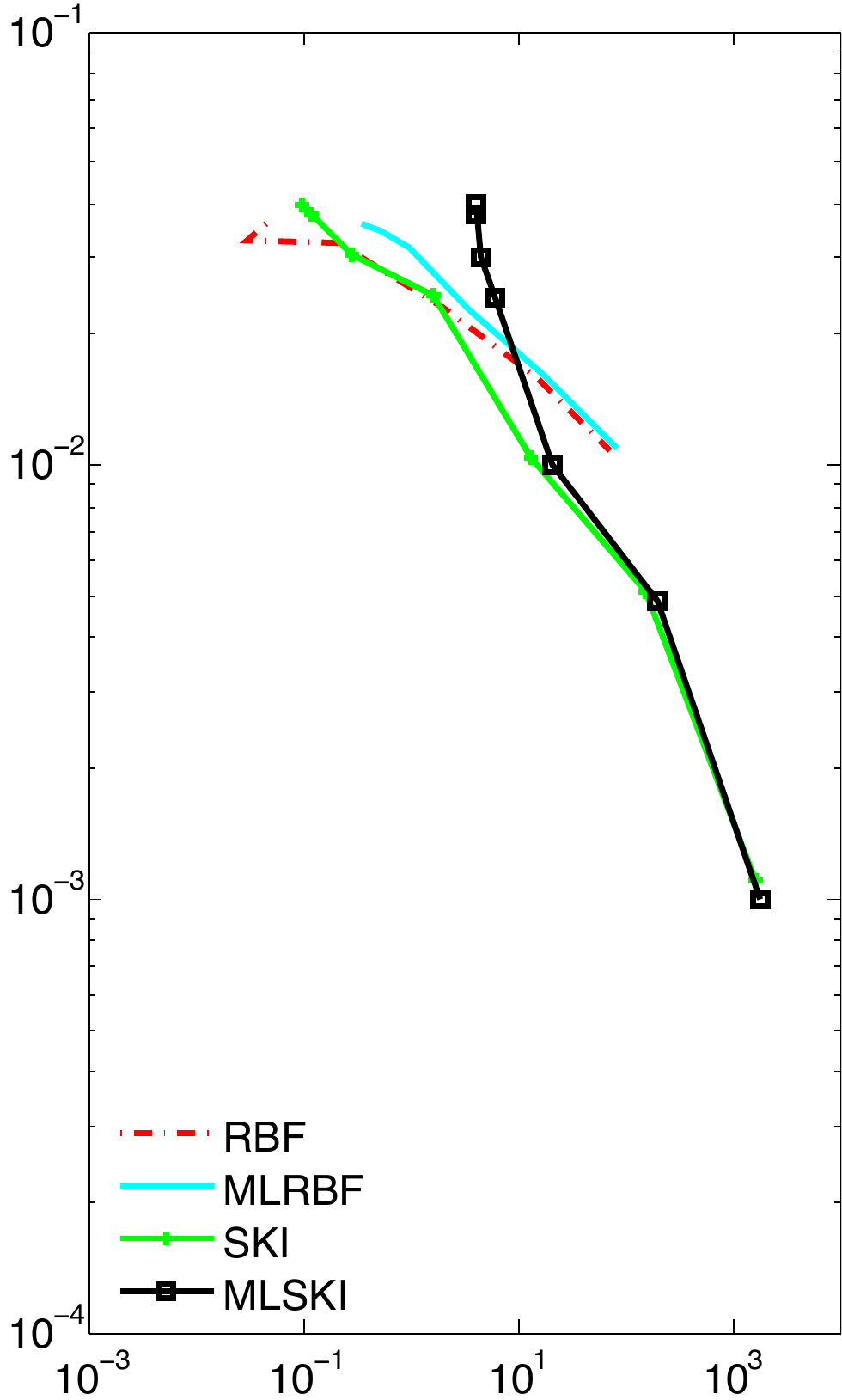}
\begin{picture}(0,0)
\put(-90,-10){Time} \put(-170,65){\rotatebox{90}{RMS-error}}
\end{picture}
}
\end{center}
\caption{Convergence of Gaussians for $u_{F4D}$. Error evaluated at
$194,481$ Halton points.} \label{F4D}
\end{figure}

\begin{figure}
\begin{center}
\subfigure{
\includegraphics[width=5.5cm,height=6cm]{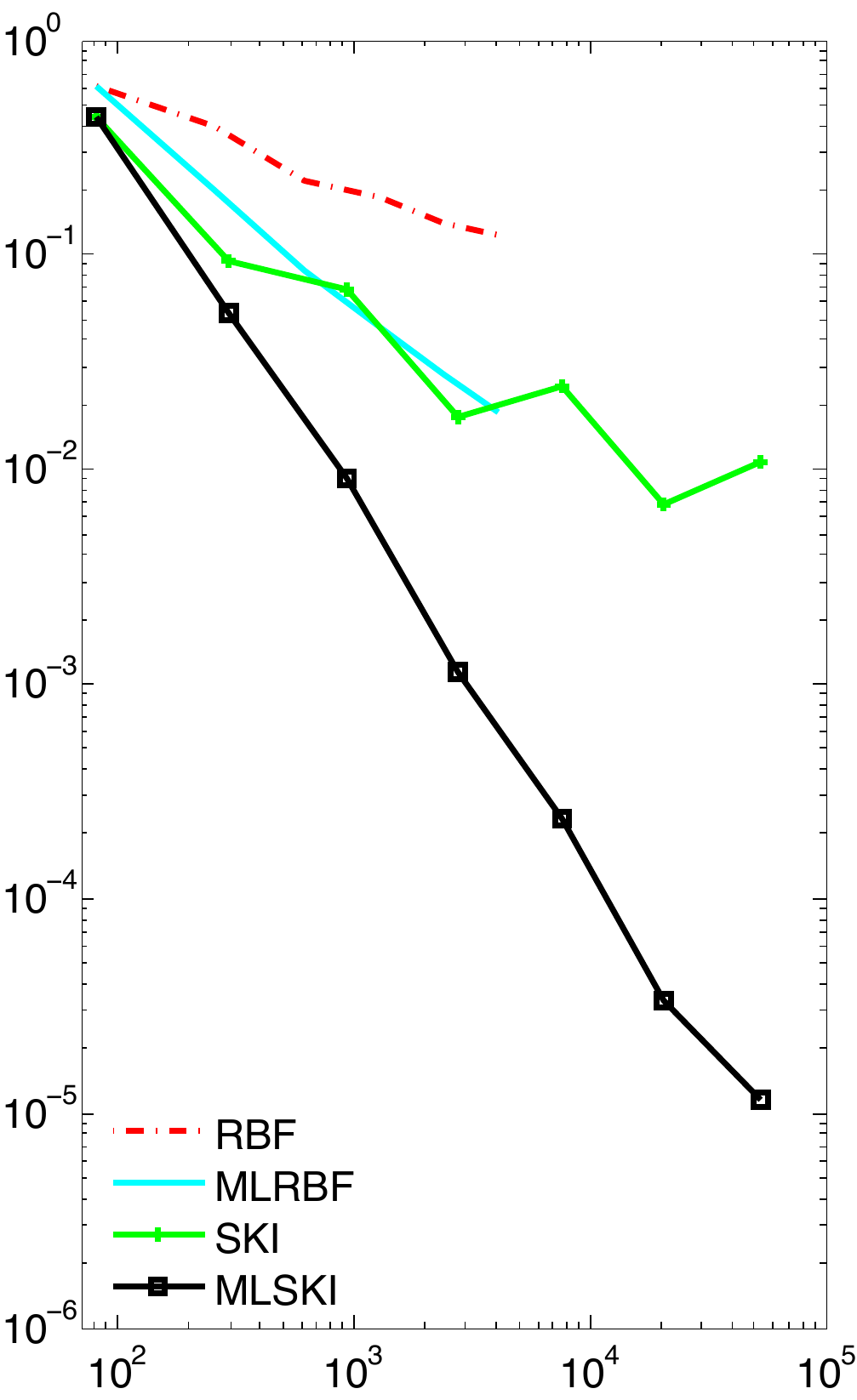}
\begin{picture}(0,0)
\put(-80,-10){$N$} \put(-170,65){\rotatebox{90}{RMS-error}}
\end{picture}
}
\subfigure{
\includegraphics[width=5.5cm,height=6cm]{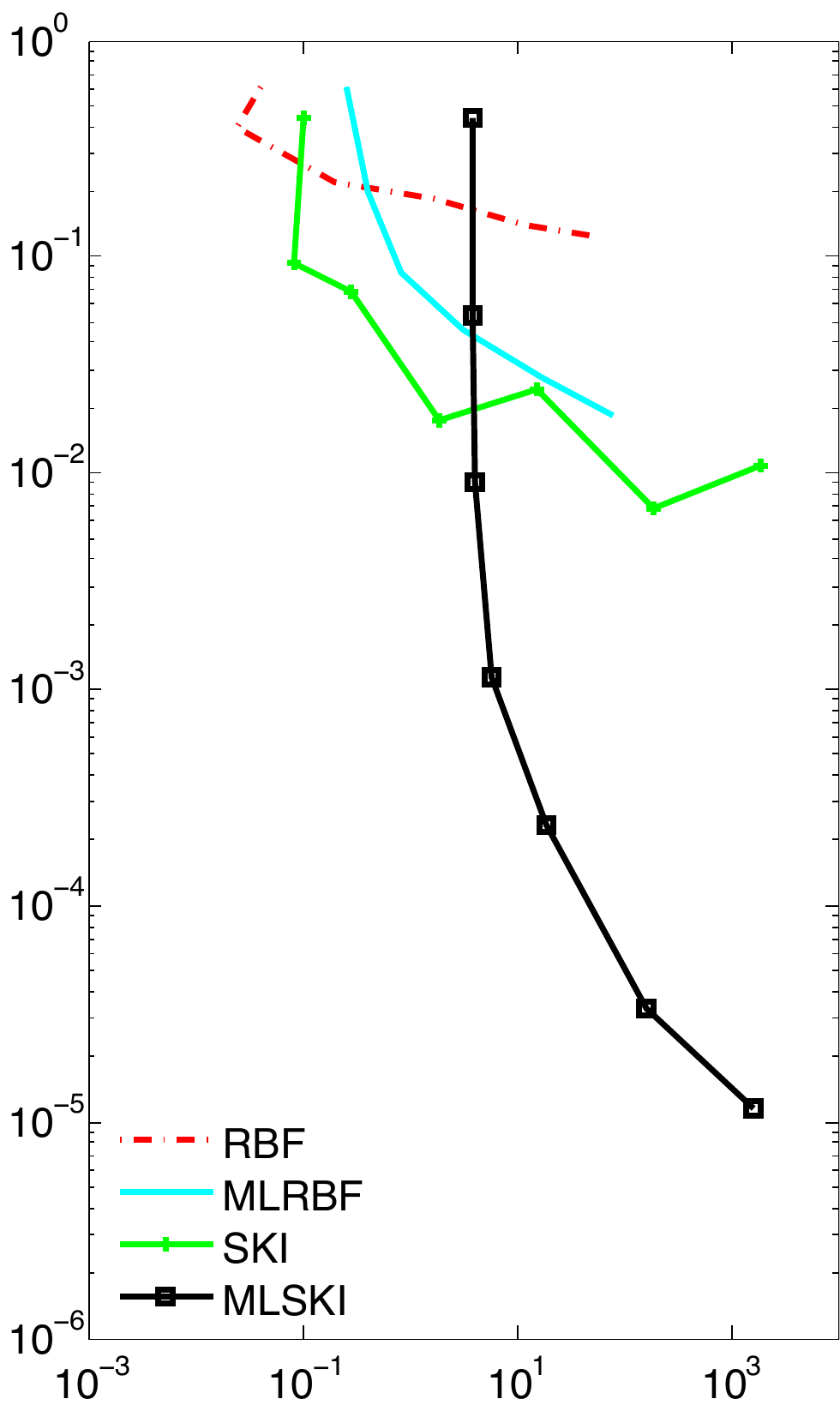}
\begin{picture}(0,0)
\put(-90,-10){Time} \put(-170,65){\rotatebox{90}{RMS-error}}
\end{picture}
}
\end{center}
\caption{Convergence of Gaussians for $u_{R4D}$. Error evaluated at
$194,481$ Halton points.} \label{R4D}
\end{figure}

We compare the four interpolation methods on a four-variate version of the Franke's function $u_{F4D}:[0,1]^4\to\mathbb{R}$, with
 \begin{equation*}
\begin{aligned}
   u_{F4D}(x_{1},\dots,x_{4})  := &\; \frac{3}{4}e^{(-(9x_{1}-2)^{2}-(9x_{2}-2)^{2} -(9x_{3}-2)^{2})/4 -(9x_{4}-2)^{2})/8}\\
     &  +\frac{3}{4}e^{-((9x_{1}+1)^{2})/49-((9x_{2}+1)^{2})/10-((9x_{3}+1)^{2})/29-((9x_{4}+1)^{2})/39}\\
      & +\frac{1}{2}e^{-((9x_{1}-7)^{2})/4-(9x_{2}-3)^{2}-((9x_{3}-5)^{2})/2-((9x_{4}-5)^{2})/4}\\
       & -\frac{1}{5}e^{-((9x_{1}-4)^{2})/4-(9x_{2}-7)^{2}-((9x_{3}-5)^{2})-((9x_{4}-5)^{2})},
\end{aligned}
 \end{equation*}
as well as the four-variate version $u_{R4D}$ of $u_{R3D}$, which is defined completely analogously. The convergence history, given in Figures \ref{F4D} and \ref{R4D}, respectively, indicates similar behaviour to the three-dimensional case. The choice in the shape parameter is as in \eqref{shape_choice} with $K=3$. We remark on the favourable complexity of the MLSKI algorithm compared to RBF/MLRBF as $d$ grows. Indeed, the MLSKI algorithm is able to compute highly accurate interpolants, as large $N$ can be achieved using moderate computational time.

\subsection{Experiment 4} It is easy to see that $d$-dimensional Gaussian kernels are tensor-products of one-dimensional ones. Hence, it is interesting to also consider SKI and MLSKI with non-Gaussian kernels. 

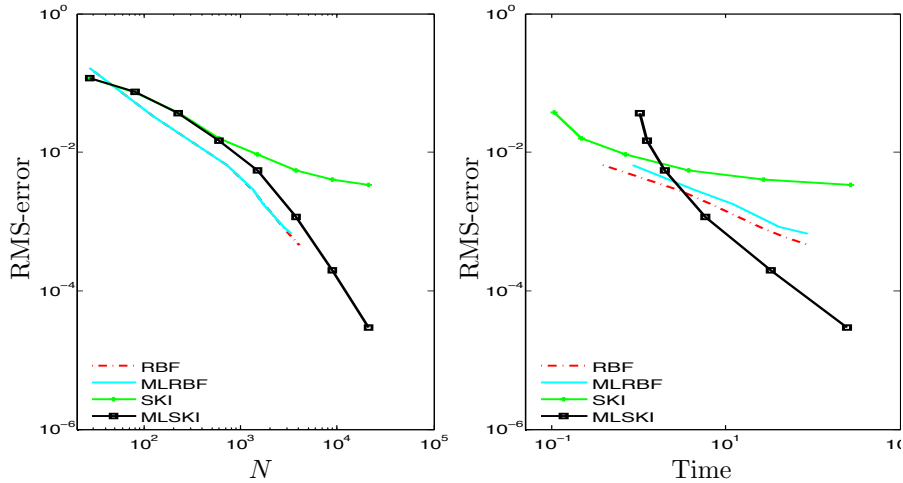
\begin{figure}
\begin{center}
\subfigure{
\includegraphics[width=5.5cm,height=6cm]{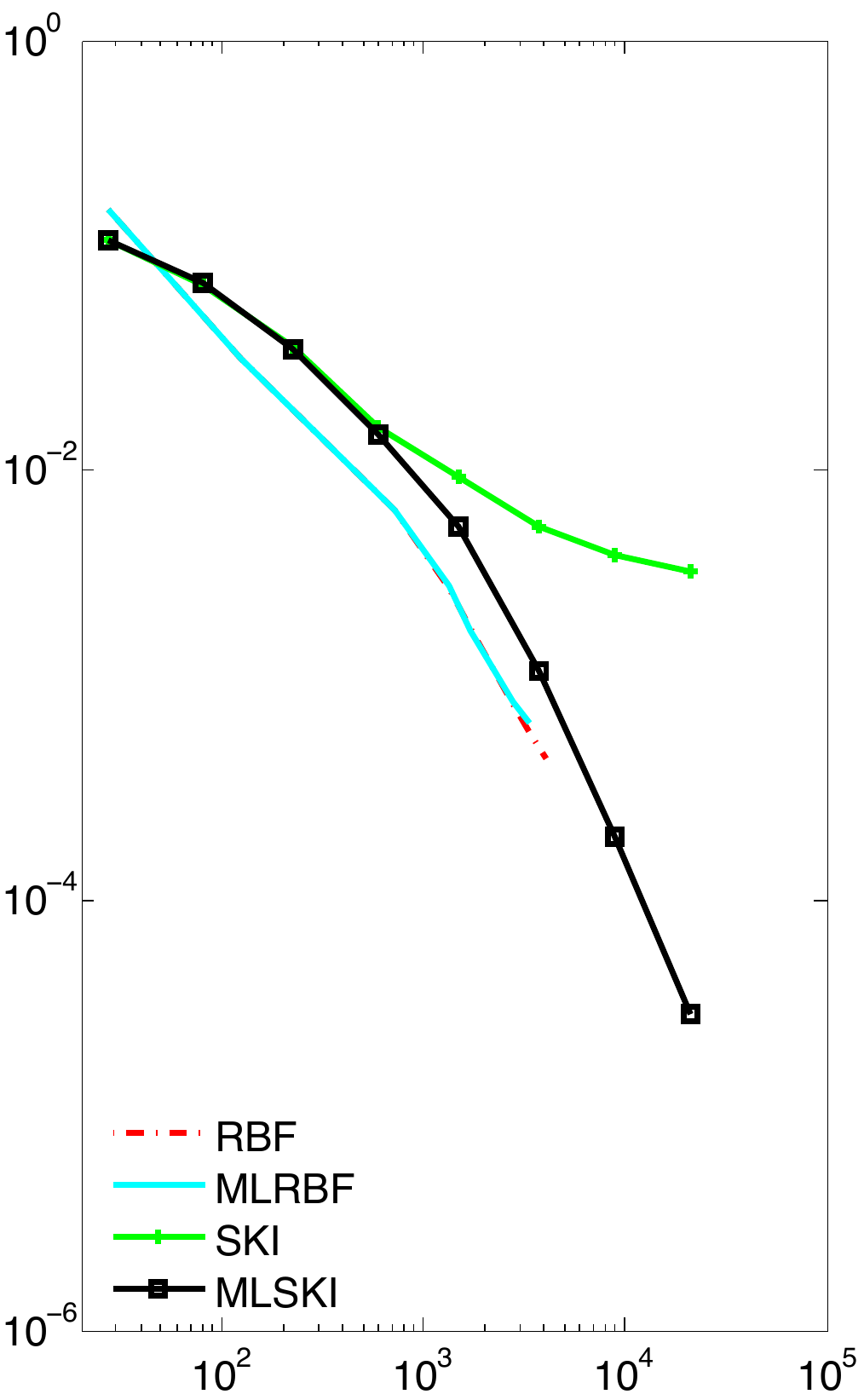}
\begin{picture}(0,0)
\put(-80,-10){$N$} \put(-170,65){\rotatebox{90}{RMS-error}}
\end{picture}
}
\subfigure{
\includegraphics[width=5.5cm,height=6cm]{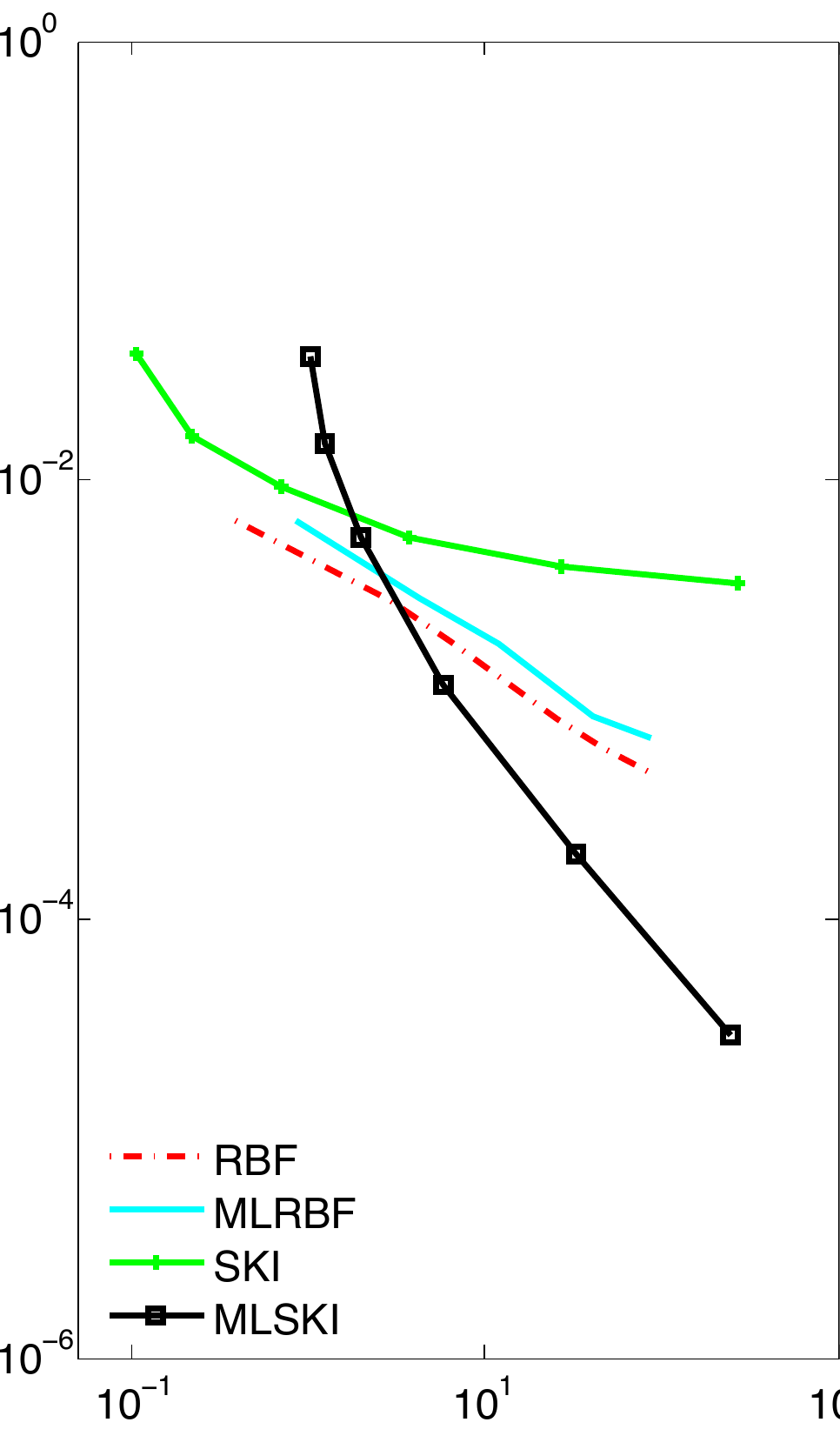}
\begin{picture}(0,0)
\put(-90,-10){Time} \put(-170,65){\rotatebox{90}{RMS-error}}
\end{picture}
}
\end{center}
\caption{Convergence of $\phi_{3,2}$ for $u_{F3D}$. Error evaluated at
$125,000$ Halton points.} \label{F3D_WE32}
\end{figure}

To this end, we consider the problem of interpolation of $u_{F3D}$ using Wendland's compactly the supported kernel $\phi_{3,2}(r):=(1-cr)^6_+(35c^2r^2+18cr+3)$.  In Figure \ref{F3D_WE32} the convergence history of the various interpolation methods based on the compactly supported kernel $\phi_{3,2}$ are given. The choice in the shape parameter is as in \eqref{shape_choice} with $K=3$, resulting in safe conditions numbers for all methods. Interestingly, the SKI method seems to converge very slowly, while the MLSKI algorithm appears to perform well.

\begin{figure}
\begin{center}
\subfigure{
\includegraphics[width=5.5cm,height=6cm]{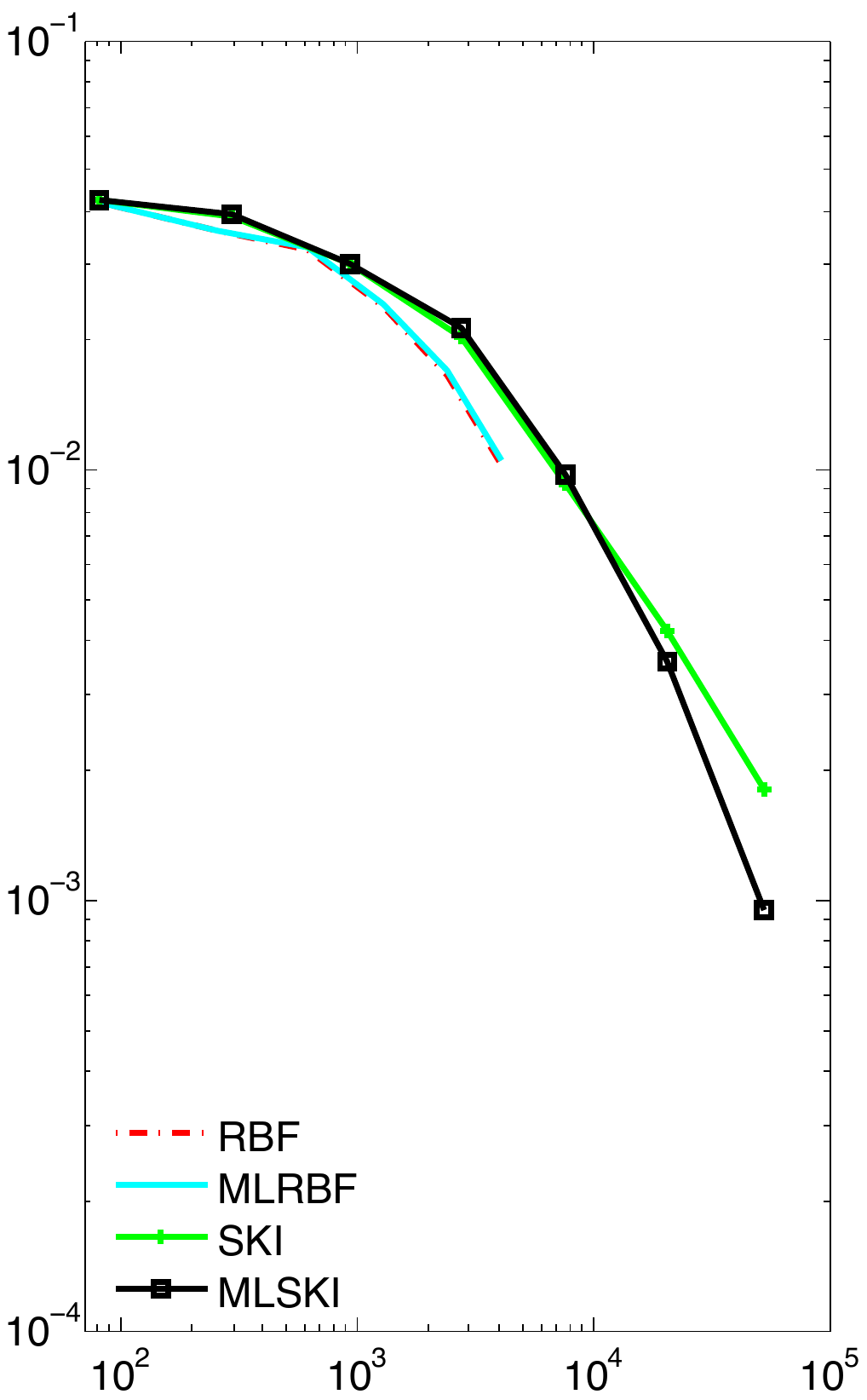}
\begin{picture}(0,0)
\put(-80,-10){$N$} \put(-170,65){\rotatebox{90}{RMS-error}}
\end{picture}
}
\subfigure{
\includegraphics[width=5.5cm,height=6cm]{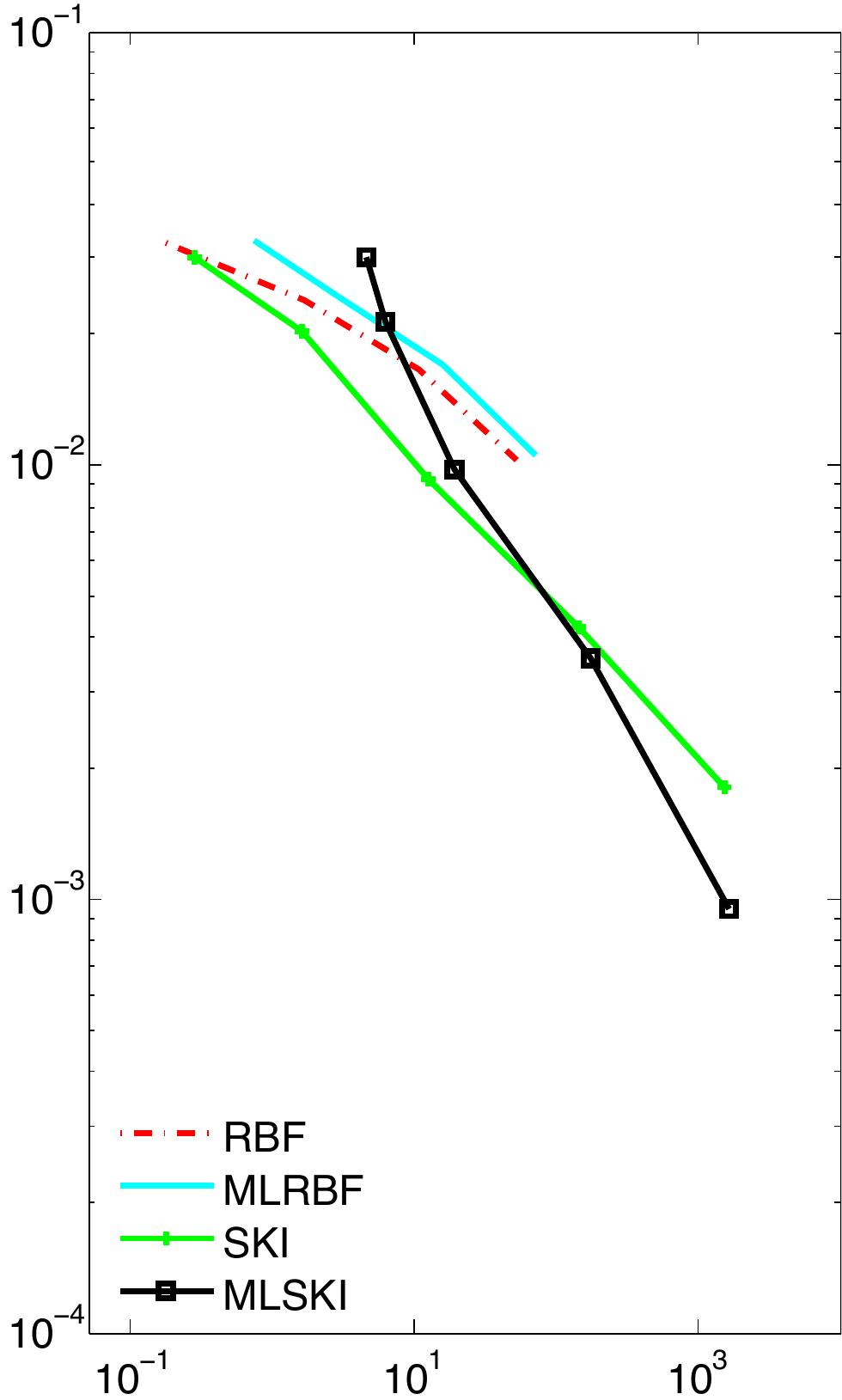}
\begin{picture}(0,0)
\put(-90,-10){Time} \put(-170,65){\rotatebox{90}{RMS-error}}
\end{picture}
}
\end{center}
\caption{Convergence of $\phi_{IMQ}$ for $u_{F4D}$. Error evaluated at
$194,481$ Halton points.} \label{F4D_IMQ}
\end{figure}
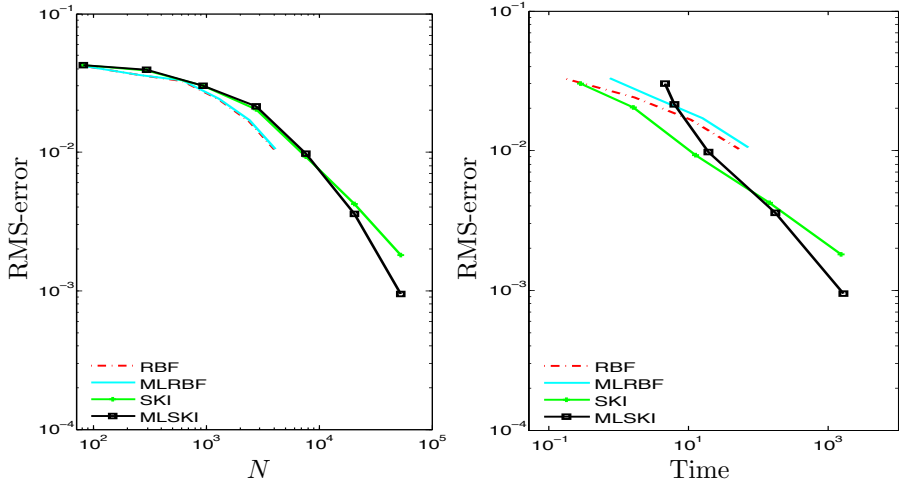

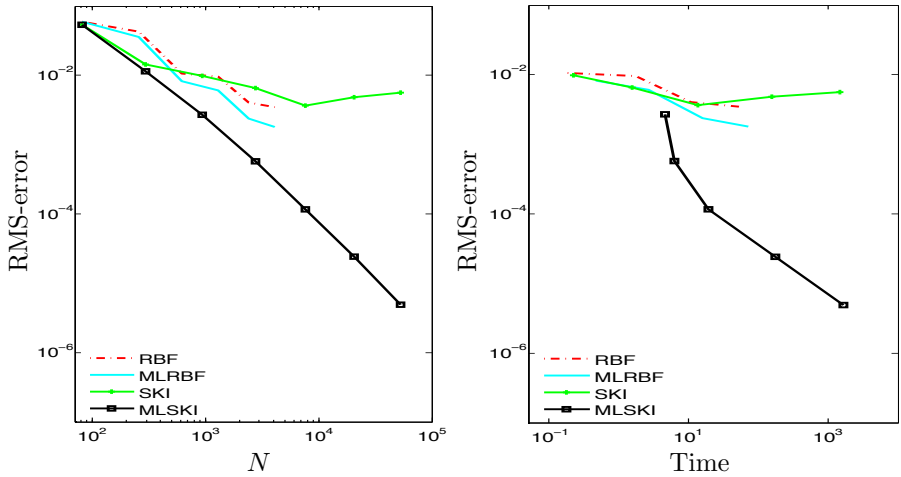
\begin{figure}
\begin{center}
\subfigure{
\includegraphics[width=5.5cm,height=6cm]{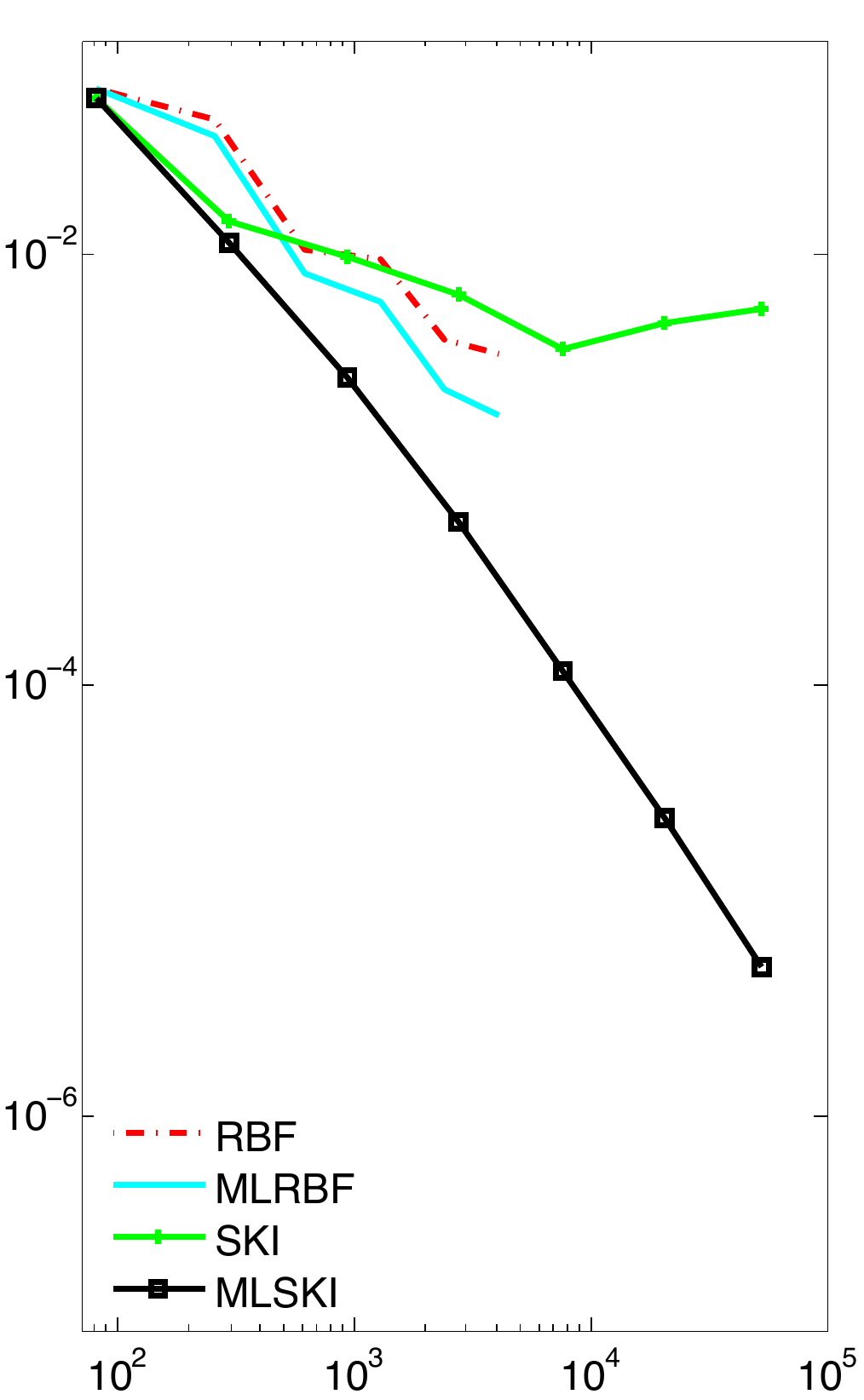}
\begin{picture}(0,0)
\put(-80,-10){$N$} \put(-170,65){\rotatebox{90}{RMS-error}}
\end{picture}
}
\subfigure{
\includegraphics[width=5.5cm,height=6cm]{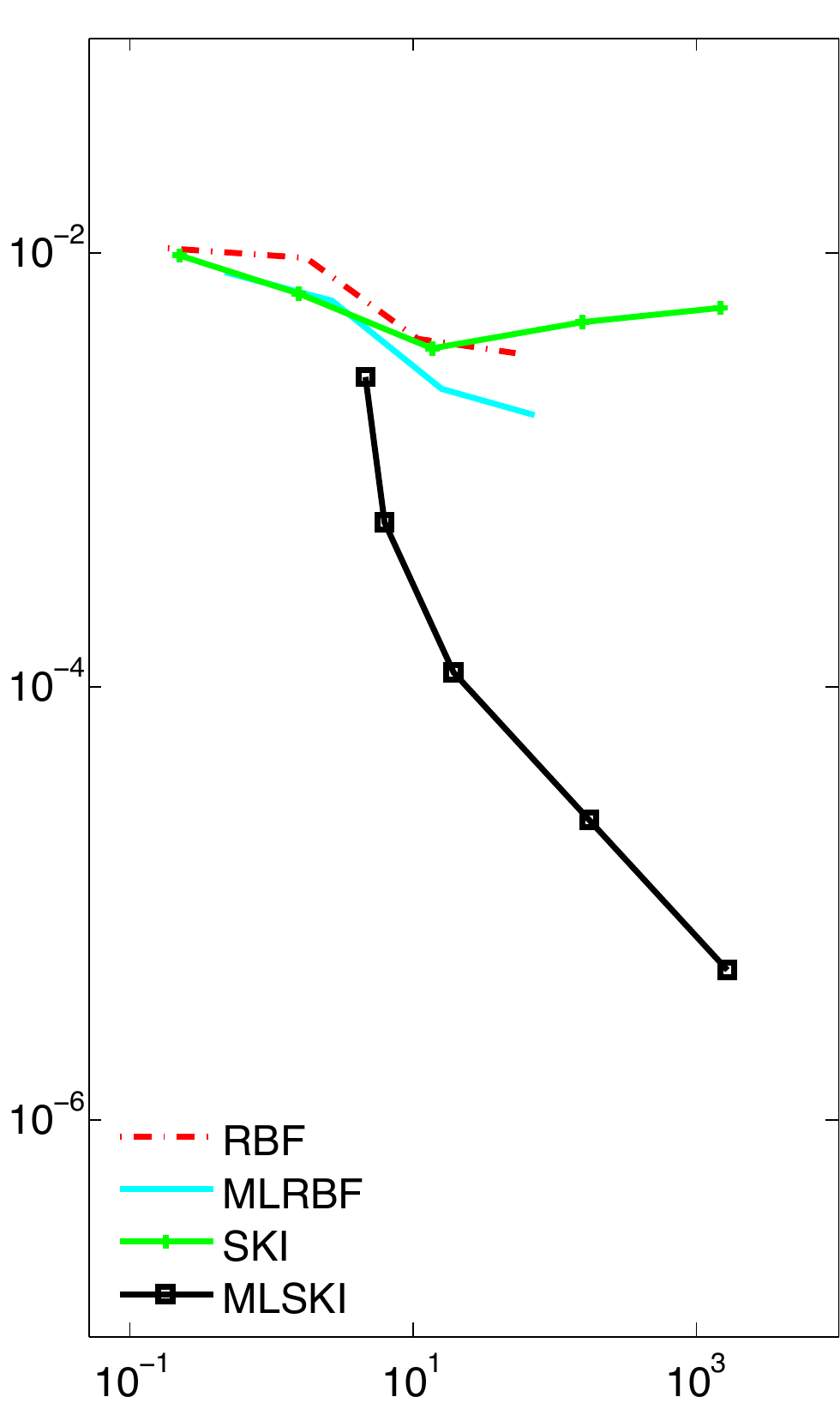}
\begin{picture}(0,0)
\put(-90,-10){Time} \put(-170,65){\rotatebox{90}{RMS-error}}
\end{picture}
}
\end{center}
\caption{Convergence of $\phi_{IMQ}$ for $u_{quad}$, for $d=4$. Error evaluated at
$194,481$ Halton points.} \label{QUAD_IMQ}
\end{figure}

Finally, we consider  the problem of interpolation of $u_{F4D}$  and  $u_{quad}$ using the inverse multiquadric kernel $\phi_{IMQ}(r):=(1+c^2r^2)^{-1/2}$ with the above choice of the shape parameter. The results are given in Figures \ref{F4D_IMQ} and \ref{QUAD_IMQ}, respectively. Again, for the case of $u_{quad}$ the SKI methods seems to be performing poorly, compared to the fast convergence of MLSKI.

\section{Concluding remarks}\label{sec7}
A multilevel kernel-based interpolation method, 
suitable for moderately high-dimensional interpolation problems on carefully structured grids has been proposed. The key idea is the use of hierarchical decomposition of the data sites with anisotropic radial basis functions used on each level, solving a number of smaller independent interpolation problems. MLSKI appears to be generally superior over classical radial basis function methods in terms of complexity, run time and convergence, at least for large data sets when $d=3,4$. It is expected that good convergence can be obtained for $d= 5$ also.

Currently, the choice of data sets (sparse grids) is highly structured. Some preliminary numerical experiments for SKI on mildly perturbed sparse grids, presented in \cite{Fazli_thesis}, indicate that the SKI method converges with the same rate, albeit with a somewhat larger constant. Perhaps a more robust methodology for extending the applicability of SKI/MLSKI methods to scattered data is the pre-computation of the values on the corresponding sparse grid data sites via local interpolation. The extension of the SKI/MLSKI method to more general geometries could possibly be handled either by introducing fictitious gridded data sites with suitable data values, taking into account the nature of the data, or by conformally mapping the computational domain \cite{MR2671290}.

The discussion in this work has been confined to strictly positive definite kernels. The implementation of MLSKI with conditionally positive definite kernels is subject to ongoing work; some preliminary results can be found in \cite{Fazli_thesis}.

\def\cprime{$'$}

\end{document}